\documentclass[11pt]{article}
\usepackage[dvips]{graphicx}
\usepackage{latexsym}
\usepackage{amssymb}

\newcommand{\C}{\mathbb{C}}
\newcommand{\CP}{\mathbb{CP}}

\newcommand{\R}{\mathbb{R}}

\renewcommand{\d}{\mathrm{d}}

\newcommand{\koniec}{\hfill $\Box $}
\def\be{\begin{equation}}
\def\ee{\end{equation}}

\def\t{\tilde}

\def\Sm{\Sigma}
\def\Smt{\widetilde {\Sm}}

\def\Th{\Theta}

\def\O{\cal O}
\def\om{\omega}

\def\ov{\overline}

\def\tw{\tilde w}
\def\tz{\tilde z}
\def\lt{\tilde{\lambda}}
\def\p{\partial}

\def\ov{\overline}

\def\Kt{\tilde{K}}
\def\K{\kappa}

\newcommand{\hook}{{\setlength{\unitlength}{11pt}   
                   \begin{picture}(.833,.8)
                   \put(.15,.08){\line(1,0){.35}}
                   \put(.5,.08){\line(0,1){.5}}
                   \end{picture}}}
\def\a{\alpha}
\def\l{\lambda}

\def\O{{\cal O}}

\def\Kt{\widetilde{K}}
\newtheorem{theo}{Theorem}[section] 
\newtheorem{prop}[theo]{Proposition}  
\newtheorem{lemma}[theo]{Lemma}

\newtheorem{col}[theo]{Corollary}

\begin{document}
\title{Einstein--Weyl geometry, the dKP equation and twistor theory}

\author{Maciej  Dunajski\thanks{email: dunajski@maths.ox.ac.uk},\,\,
Lionel J. Mason,\,\,  Paul Tod\\ The Mathematical Institute,
24-29 St Giles, Oxford OX1 3LB, UK}  

\date{} 
\maketitle
\abstract { It is shown that Einstein--Weyl (EW) equations in 2+1 dimensions contain the
dispersionless Kadomtsev--Petviashvili (dKP) equation  as a special
case:
If an EW structure admits a constant weighted vector then it is
locally given by
$h=\d y^2-4\d x\d t-4u\d t^2, \nu=-4u_x\d t$, where
$u=u(x, y, t)$ satisfies 
the dKP equation
$(u_t-uu_x)_x=u_{yy}$. 

Linearised solutions to the dKP equation are shown to give rise to
four-dimensional anti-self-dual conformal structures with symmetries.
All four-dimensional hyper-K\"ahler metrics in signature $(++--)$
for which the self-dual part of the derivative of a Killing vector is
null arise by this construction.

Two new  classes of examples of EW
metrics which depend on one arbitrary function of one
variable are given, and characterised.  

A Lax representation of the EW condition is found and used
to show that all EW spaces arise as symmetry reductions of 
hyper-Hermitian metrics in four dimensions.

The  EW equations are reformulated 
in terms of a simple and closed two-form on the $\CP^1$-bundle over a 
Weyl space. 

It is proved that complex solutions to the dKP equations, modulo a
certain coordinate freedom, are in a one-to-one correspondence with
minitwistor spaces (two-dimensional complex manifolds ${\cal Z}$
containing a rational curve with normal bundle $\O(2)$) that admit a
section of $\K^{-1/4}$, where $\K$ is the canonical bundle of ${\cal
Z}$.  Real solutions are obtained if the minitwistor space also admits
an anti-holomorphic involution with fixed points together with a
rational curve and section of $\K^{-1/4}$ that are invariant under the
involution.}

\section{Three-dimensional Einstein--Weyl spaces}
The aim of this paper is to study the Einstein--Weyl (EW) equations in
relation to integrable systems, and in particular the dispersionless
Kadomtsev--Petviashvili equation.

We begin by collecting various definitions and formulae
concerning three-dimensional Einstein--Weyl spaces (see \cite{PT93}
for a fuller account).  In section 2 we construct and
characterise a class of new EW structures in 2+1 dimensions out of
solutions to the dKP equation.  We then show that the dKP
solutions give rise to hyper-K\"ahler metrics in four dimensions. We
abuse terminology and call hyper-K\"ahler (hyper-complex,
hyper-Hermitian) metrics which in signature $(++--)$ should be
referred to as pseudo-hyper-K\"ahler (pseudo-hyper-complex,
pseudo-hyper-Hermitian). A null vector field 
(with conformal weight) will play a central role in our discussion 
so most of our constructions only make sense for
Einsetin-Weyl spaces with Lorentzian signature, or complex holomorphic
EW spaces (i.e.\ the complexification of real analytic EW spaces) and
for the most part we work with the latter and restrict to a real slice
when reality conditions play a role.

In section 3 we construct some new examples of EW structures.
We obtain all solutions of the dKP equation with the property
that the associated EW space admits a family of divergence-free, shear-free 
geodesic congruences. These solutions give rise to new EW metrics
depending on one arbitrary function of one variable.

In section 4 a Lax representation of the
general EW equations is given, together with a reformulation of 
the EW equations in terms of a closed and simple two-form on the bundle
of spinors. 
A full twistor characterisation of dKP Einstein--Weyl structures and
the corresponding hyper-K\"ahler metrics will be given in section 5.
In section 6 we summarise our present knowledge of conformal
reductions of four-dimensional hyper-K\"ahler metrics in split signature.
In the Appendix we show how to obtain the dKP equation as a
reduction of Pleba{\'n}ski's second heavenly equation \cite{Pl75}.
Parts of this work appeared in the DPhil thesis of one of the authors 
(MD) \cite{DPhil}. 
 
 Let ${\cal W}$ be a
$3$-dimensional complex manifold (one can also define Weyl spaces in
arbitrary dimension) with a torsion-free connection $D$ and a
conformal metric $[h]$.  We shall call ${\cal W}$ a Weyl space if the
null geodesics of $[h]$ are also geodesics for $D$.  This condition is
equivalent to \be
\label{ew1}
D_ih_{jk}=\nu_ih_{jk}
\ee
for some one form $\nu$. Here $h_{jk}$ is a representative metric in
the conformal class. The indices 
$i, j, k, ...$ go from $1$ to $3$.
If we change this representative by
$h\longrightarrow \phi^2 h$, then $\nu\longrightarrow
\nu+2\d\ln{\phi}$.  
The one-form $\nu$ `measures' the difference between $D$ and
the Levi-Civita connection $\nabla$ of $h$:
\be
\label{difference}
D_iV^j=\nabla_iV^j-\frac{1}{2}\Big(\delta_{i}^j\nu_k +\delta_{k}^j\nu_i    
-h_{ik}\nu^j\Big)V^k.
\ee
The Ricci tensor $W_{ij}$ and scalar $W$ 
of $D$ are related to the Ricci tensor
$R_{ij}$ and scalar $R$ of $\nabla$ by
\begin{eqnarray*}
W_{ij}&=&R_{ij}+\nabla_i\nu_j-\frac{1}{2}\nabla_j\nu_i
+\frac{1}{4}\nu_i\nu_j+h_{ij}\Big(-\frac{1}{4}\nu_k\nu^k+\frac{1}{2}
\nabla_k\nu^k\Big),\\
W:&=&h^{ij}W_{ij}=R+2\nabla^k\nu_k-\frac{1}{2}\nu^k\nu_k.
\end{eqnarray*}
A tensor object $T$ which transforms as $
T\longrightarrow \phi^m T$ when
$ h_{ij}\longrightarrow \phi^2 h_{ij}$
is said to be conformally invariant of weight $m$.
The Ricci scalar $W$, and the Ricci tensor $W_{ij}$ have weights 
$-2$ and $0$ respectively.

Let $\beta$ be a $p$-form of weight $m$. The covariant exterior 
derivative
\[
\widetilde{D}\beta:=\d\beta-\frac{m}{2}\nu\wedge\beta
\]
is a well-defined $p+1$-form of weight $m$. The formula for a
covariant weighted derivative of a vector of weight $m$ is
\be
\label{wcowder}
\widetilde{D}_iV^j= \nabla_iV^j-\frac{1}{2}\delta_i^j\nu_kV^k
-\frac{m+1}{2}\nu_iV^j+\frac{1}{2}\nu ^jV_i.
\ee
We say that a vector $K$ is a symmetry of a Weyl structure if 
it preserves the conformal structure $[h]$, the  Weyl connection, and 
the compatibility (\ref{ew1}) between those two. These conditions imply
\be
\label{Weyl_Killing}
{\cal L}_Kh=\psi h,\qquad {\cal L}_K\nu=\d\psi,
\ee
where $(h, \nu)$ is a Weyl structure, and ${\cal L}_K$ is the Lie
derivative along $K$. 

The conformally invariant Einstein--Weyl (EW) condition on
$({\cal W}, h, \nu)$ is
\[
W_{(ij)}=\frac{1}{3}W h_{ij}.
\]
If the above equation is satisfied and $ \nu$ is a gradient, 
then $h$ is conformal to a metric with
constant curvature. 

In terms of the Riemannian data the Einstein--Weyl equations are
\be
\label{ew2}
\chi_{ij}:=R_{ij}+\frac{1}{2}\nabla_{(i}\nu_{j)}+\frac{1}{4}\nu_i\nu_j
-\frac{1}{3}\Big(R+\frac{1}{2}\nabla^k\nu_k+\frac{1}{4}\nu^k\nu_k\Big)h_{ij}=0.
\ee
Here $\chi_{ij}$ is a conformally invariant tensor
(the trace-free part
of the Ricci tensor of the Weyl connection).
Weyl spaces which satisfy (\ref{ew2}) will be called Einstein--Weyl
(or EW) spaces.

In three dimensions  the general solution of (\ref{ew1})-(\ref{ew2}) depends on four arbitrary functions
of two variables \cite{C43}. 
The equations of the Weyl geodesics are
\[
\frac{\d}{\d s}\frac{\p {\cal L}}{\p{\dot x}^i}-\frac{\p {\cal L}}{\p x^i}
=F_i(x^j, \dot x^j)
\]
where ${\cal L}= (1/2)h_{ij}{\dot x}^i{\dot x}^j$ and 
$F_i={\dot x}_i({\dot x}^j\nu_j)-(1/2)\nu_i({\dot x}^j{\dot x}_j)$. 
Here $\dot{}=\d/\d s$ stands
for the derivative with respect to a parameter $s$.
It is evident that for null $\dot{x}^i$ the geodesics 
coincide with the null geodesics
for $[h]$.

\section{Einstein--Weyl structures from the dKP equation}
In this section we shall construct
Einstein--Weyl structures out of solutions to the dKP equation. 
In subsection 2.1 we shall 
find a class of hyper-K\"ahler metrics in four dimensions which
reduce to dKP EW metrics.

The full  Kadomtsev--Petviashvili equation for  $U:=U(X^i), X^i=(X, Y, T)$ 
\be
\label{KP}
(U_T-UU_X-(1/12)U_{XXX})_X=U_{YY}
\ee
arises as a compatibility condition for the linear system $L_0\Psi=L_1\Psi=0$,
where $\Psi=\Psi(X, Y, T)$ and
\[
L_0=\p_Y-(1/2)\p_X^2-U,\qquad L_1=\p_T-(1/3)\p^3_X -U\p_X-W,
\]
for some $W=W(X, Y, T)$.  
To take a dispersionless limit of (\ref{KP}) \cite{G85} 
introduce the slow coordinates $x^i:=\epsilon X^i$ (note that our notation for
`slow' and `fast' coordinates is different from the usual one), and define
$u(x^i):=U(X^i), w(x^i):=W(X^i)$.
The linear system is replaced by
\be
\label{HamJac}
S_y=(1/2)S_x^2+u,\qquad S_t=(1/3)S_x^3+uS_x+w.
\ee
Here $S:=S(x^i)$ is the action defined by $\Psi(X^i)=\exp{[\epsilon^{-1}S(x^i)]}$,
and higher order terms in $\epsilon$ have been neglected.
 Formulae (\ref{HamJac}) can be treated as a pair of Hamilton--Jacobi equations
$S_{t_A}+H_A(S_x, x, t_A)=0$, with $t_A =(y,t)$ and $H_A=(H_2, H_3)$ where
\[
H_2:=\frac{{\tilde\lambda}^2}{2}+u,\;\;\;H_3:=\frac{{\tilde\lambda}^3}{3}+
\tilde{\lambda}{u}+w
\]
for  $u=u(x, y, t)$ and $w=w(x, y, t)$.

Now $x^i$ and  $\p S/\p x^i=(\lt, H_2, H_3)$ form a set of canonically
conjugate variables
on an `extended phase-space', with the symplectic form
\be
\label{dKPtwoform}
\Pi=\d x^i\wedge\d\frac{\p S}{\p x^i}=\d x\wedge\d \lt+\d y\wedge\d
H_2+ \d t\wedge \d H_3.  
\ee 
This two-form is closed by definition. It is also simple iff $u$ and
$w$ satisfy 
\[
w_x=u_y,\qquad u_t-uu_x=w_y.
\]
Eliminating $w$ yields the 
dKP equation
\be
\label{dKP}
(u_t-uu_x)_x=u_{yy}.
\ee
The simplicity of $\Pi$ implies 
$
[\p_y+ X_{H_2}, \p_t+ X_{H_3}]=0
$
where $X_H:=H_x\p_{\lt}-H_{\lt}\p_x$ denotes the Hamiltonian vector field with respect to 
$\d \t{\l}\wedge\d x$, holding $t$ and $y$ constant. This gives a Lax pair for the dKP equation in terms of
Hamiltonian vector fields. To obtain a Lax pair which is linear in the
spectral parameter put
\be
\label{dKPlax}
L_{0'}:=\p_t+X_{H_3}-\lt(\p_y+ X_{H_2})
=\p_t-u\p_x-\t{\l} \p_y+u_y\p_{\t{\l}},\qquad
L_{1'}:=\p_y+X_{H_2}=\p_y-\lt \p_x+u_x\p_{\t{\l}}.
\ee
The dKP equation is equivalent to 
\[
[L_{0'}, L_{1'}]=-u_xL_{1'}.
\]
Define a triad of vectors
\[
\nabla_{1'1'}:=\p_x,\;\;\;\;\nabla_{0'1'}:=\p_y,\;\;\;\;
\nabla_{0'0'}:=\p_t-u\p_x
\]
so $L_{A'}=\pi^{B'}\nabla_{A'B'}+f_{A'}\p_{\lt}$, 
where $\pi^{A'}=(1, -\lt)$ and $f_{A'}=(u_y, u_x)$.

The next proposition shows that we can find a one form $\nu$ such
that $\nabla_{A'B'}$ is a null triad for an EW metric:
\begin{prop}
Let $u:=u(x,y,t)$ be a solution of the {\em dKP} equation {\em(\ref{dKP})}.
Then the metric and the one-form
\be
\label{EWdkp}
h=\d y^2-4\d x\d t-4u\d t^2,\qquad\nu=-4u_x\d t
\ee
give an EW structure.
\end{prop}
{\bf Proof.} Let $x^1:=t, x^2:=y, x^3:=x$. 
Five (out of six) EW equations $\chi_{ij}=0$ are satisfied identically
by ansatz (\ref{EWdkp}). The equation $\chi_{11}=0$ is equivalent to
(\ref{dKP}). We also find $W=-3u_{xx}$.\koniec

\smallskip

\noindent
{\bf Example:} Solutions which yield EW structures conformal to
Einstein metrics (i.e. those for which $\nu$ is exact) 
are of the form 
\be
\label{conformal_Einstein}
u(x,y,t)=xf_1(t)+
\frac{1}{2}\Big(\frac{\d {f_1}(t)}{\d t}-f_1(t)^2\Big)y^2+f_2(t)y+f_3(t),
\ee
where $f_1(t), f_2(t), f_3(t)$ are arbitrary functions of one
variable.

One can verify that the  vector $\p_x$ in the EW space (\ref{EWdkp})
is a covariantly constant null vector in the Weyl connection  
with weight $-1/2$.
Now we shall prove the converse, and show 
that solutions (\ref{EWdkp}) are characterised by
the existence of a constant weighted vector.
\begin{prop}
\label{todth}
If a three dimensional EW space has a constant weighted vector field $l$
then coordinates can be
chosen to put the EW metric and 1-form in the form  {\em (\ref{EWdkp})}.
\end{prop}
We shall need  following lemma:
\begin{lemma}
\label{lemma_weight}
Let $l$ be a constant weighted vector on a three-dimensional EW space.  
Then either the EW space 
is flat or $l$ is null  
{\em(}so on a real slice the signature is $(+--)${\em)} and has weight $-1/2$.
\end{lemma}
{\bf Proof.} 
Assume that $(h, \nu)$ is a complex EW structure (we shall specify the
reality conditions later in the proof).
Commuting the Weyl derivatives yields
\[
[D_i, D_j]l^k=\frac{m}{2}(D_i\nu_j-D_j\nu_i)l^k={W^k}_{mij}l^m,
\]
where ${W^k}_{mij}$ is the curvature of the Weyl connection, and $m$
is the weight of $l^k$.
It can be decomposed as 
\be
\label{tt1}
{W^k}_{mij}=-{\varepsilon_{ij}}^p{\varepsilon_{m}}^{kq}S_{pq}
-\delta_m^kF_{ij},
\ee
where $F_{ij}=\nabla_{[i}\nu_{j]}$, and $S_{ij}$ is a conformally 
invariant tensor of weight $0$. If the EW equations are satisfied
$S_{ij}$ is given by
\be
\label{tt2}
S_{ij}=\frac{1}{2}F_{ij}+\frac{ W}{6}h_{ij}.
\ee
Equations (\ref{tt1}) and  (\ref{tt2}) imply
\be
\label{tt3}
(m+1)F_{ij}l^k=-\frac{1}{2}{\varepsilon_{ij}}^pl^m
{\varepsilon_{m}}^{kq}F_{pq}+
\frac{W}{6}(\delta_i^kl_j-\delta_j^kl_i).
\ee
In three dimensions any non-zero two-form $F_{ij}$ has a non-trivial
kernel, i.e. there exists a non-zero vector $L^j$ with 
$F_{ij}L^j=0$, which implies 
\be
\label{tt31}
F_{ij}=F{\varepsilon_{ijk}}L^k
\ee
for some non-zero $F$. We have to consider three cases:
\begin{itemize}
\item 
Suppose first that $L^k$ is a null vector and
contract (\ref{tt3}) with $L^j$ to find
\be
\label{tt32}
0=-\frac{1}{2}{\varepsilon_{ij}}^p{\varepsilon_{m}}^{kq}
F{\varepsilon_{pqr}}L^rl^mL^j+\frac{W}{6}(\delta_i^kl_jL^j
-L^kl_i).
\ee
Contracting this with $L_k$ yields $Wl_jL^j=0$. If $W=0$ then 
(\ref{tt32}) implies that $l^i$ and $L^i$ are proportional, so 
$l^i$ is null. If $W\neq 0$, so that $l_jL^j=0$ then 
(\ref{tt32}) reduces to 
\[
0=\frac{1}{2}FL^ql^mL_i{\varepsilon_{m}^k}_q-\frac{W}{6}l_iL^k
\]
from which again $l^i$ is null. Therefore $l^i$ and $L^i$ 
are both null and orthogonal and so (as we work in three dimensions)
they have to be proportional. Now (\ref{tt32}) forces $W=0$.
Equation  (\ref{tt3}) is now satisfied only if $m=-1/2$.
\item
If $L^i$ is not null, we can 
choose an orthogonal frame with $F_{23}=F\neq 0$ , and $F_{12}=F_{13}=0$,
and use (\ref{tt3})
to  examine components of $F_{ij}l^k$ in this frame. This yields
\begin{eqnarray}
\label{tt5}
Wl_1&=&0, \qquad Fl^1=0,\qquad \frac{1}{2}Fl^3+\frac{1}{6}Wl_2=0,\qquad
\frac{1}{2}FV^2-\frac{1}{6}Wl_3=0,\\
(m+1)Fl^1&=&0,\qquad
(m+1)Fl^2=\frac{1}{6}Wl_3=\frac{1}{2}Fl^2,\qquad
(m+1)Fl^3=-\frac{1}{6}Wl_2=\frac{1}{2}Fl^3.\nonumber
\end{eqnarray}
Therefore $l^1=0$, and (\ref{tt5}) imply
$
(m+1/2)Fl^2=0, (m+1/2)Fl^3=0. 
$
But $l^i\neq 0$, so $m=-1/2$. Equations ( \ref{tt5}) also imply  that 
$l^i$ is null.
\item If $F=0=\d\nu=0$ 
(Einstein case) choose a conformal gauge in which $\nu=0$. Now
$D_il^j=\nabla_il^j=0$ implies $R=0$. Therefore   the metric $h$ is flat  
and $l^j$ is a constant vector.
\end{itemize}
\koniec
\smallskip

\noindent
{\bf Proof of Proposition (\ref{todth}).}
Lemma  \ref{lemma_weight} and  the formula (\ref{wcowder}) with
$m=-1/2$ imply
\be
\label{weightofl}
\tilde{D}_il^j={D}_il^j+\frac{1}{4}\nu_il^j=0.
\ee
Therefore $D_il_j=(3/4)\nu_il_j$, so $\d {\bf l}=(3/4)\nu\wedge {\bf l}$
(here  ${\bf l}$ is the one form dual to $l$). 

This implies that we can rescale the metric and hence ${\bf l}$ so
that ${\bf l}=-2\d t$ for some function $t$.  We must then have $\nu=b\d
t$ for some function $b$.  Choose coordinates $x$ and $y$ so that
$l(y)=0$ and $l(x)=1$ and $(x,y,t)$ is a coordinate system.
At this point we have
\[
h=F\d y^2+G\d y\d t-4\d x\d t -4u\d t^2,\qquad\nu=b\d t,
\]
where $F, G, b$ and $u$ are functions of $x, y, t$.
The formulae (\ref{difference}) and (\ref{weightofl})
imply $\nabla_il_j=(1/4)\nu_il_j-(1/2)\nu_jl_i$. Symmetrising this
expression yields
$\nabla_{(i}l_{j)}=-(1/4)\nu_{i}l_{j}$, which implies that
$F_x=G_x=0$, and $4u_x=-b$.
We are still free to change $x\rightarrow x+P(y, t)$, 
which gives
\[
h=F\d y^2+G\d y\d t-4(\d x+P_y\d y+P_t\d t)\d t-4u\d t^2
\qquad\nu=-4u_x\d t.
\]
We can find $K$ such that $\d\hat{y}:=\sqrt{F}\d y+K\d t$ is exact, and
eliminate the $\d \hat{y}\d t$ term in the metric by choosing
$4P_y=-2K+G/\sqrt{F}$. This (after redefining $u$ by adding to it a function
of $(\hat{y},t)$ so that $\nu$ remains unchanged) yields
the EW structure (\ref{EWdkp}).
\koniec

\medskip

\noindent
{\bf Remark:} The above coordinate conditions fix the coordinates and
$u$ only
up to the freedom $(x,y,t)
\mapsto (\tilde{x},\tilde{y},\tilde{t})$, $u(x,y,t)\mapsto
\tilde{u}(\tilde{x},\tilde{y},\tilde{t})$ 
where
\begin{eqnarray}\label{coordtrans}
(x,y,t)&=&( \tilde{x}-f' \tilde{y}- g ,\tilde{y}-2f, \tilde{t})\, ,
\nonumber \\
\tilde{u}(\tilde{x},\tilde{y},\tilde{t})&=&u( \tilde{x}-f' \tilde{y}- g
,\tilde{y}-2f, \tilde{t}) -\tilde{y} f''-f'^2-g'\, .
\end{eqnarray}
where $f$ and $g$ are arbitrary functions of $t$ and $\prime$ denotes
the derivative with respect to $t$.

Furthermore the conformal scale is only fixed up to arbitrary
functions of $t$, $h\mapsto\tilde{h}=\Omega^2 h$.  Such a rescaling
leads to a redefinition of $t$, $t\mapsto \tilde{t}$ given by
$t=c(\tilde{t})$ where $\Omega=c^{\prime -2/3}$ where now and in the
following $\prime$ denotes the derivative wrt $\tilde{t}$.  This leads
to the redefinitions $(x,y,t) \rightarrow
(\tilde{x},\tilde{y},\tilde{t})$, $u(x,y,t)\rightarrow
\tilde{u}(\tilde{x},\tilde{y},\tilde{t})$ given by
\begin{eqnarray}\label{conftrans}
(x,y,t)&=&(c^{\prime
1/3}\tilde{x}+\frac{
c''}{6c^{\prime 2/3}}\tilde{y}^2, c^{\prime 2/3}
\tilde{y},c(\tilde{t})) \, ,\nonumber \\
\tilde{u}(\tilde{x},\tilde{y},\tilde{t}) &=&c^{\prime 2/3}u(c^{\prime
1/3}\tilde{x}+\frac{c''}{6c^{\prime 2/3}}\tilde{y}^2, c^{\prime 2/3}
\tilde{y},c)  +\frac{ c''\tilde{x}}{3c'}
+\frac{\tilde{y}^2}{18}\left (\frac{3c^{\prime\prime\prime}}{c'}- 
4\left(\frac{c''}{c'}\right)^2\right) \, .
\end{eqnarray}

From the point of view of the Einstein-Weyl spaces, the
transformations above are equivalences, however from the point of view of
the dKP equations, they map one solution of the dKP equations to
another allowing one to deduce solutions depending on 3 functions of
one variable from a given solution:
\begin{col}\label{solntrans}
Let $u(x,y,t)$ be a solution to the dKP equation, then
$\tilde{u}(\tilde{x},\tilde{y},\tilde{t})$ is another solution where $\tilde {u}$ is given in
terms of either of the formulae {\em(\ref{conftrans})} or
{\em(\ref{coordtrans})}.
\end{col}

\subsection{Hyper-K\"ahler structures from the dKP equation}
\label{HKsection} 
In this subsection we shall show that EW structures given by (\ref{EWdkp})
give rise to four-dimensional hyper-K{\"a}hler structures with symmetry.
We shall start by 
summarising some results about anti-self-dual (ASD) four manifolds with
Killing vectors, and the Lax representation of hyper-Hermitian four manifolds.

All three-dimensional EW spaces can be obtained as spaces of trajectories of
conformal Killing vectors in four-dimensional manifolds with ASD
conformal curvature:
\begin{prop}[\cite{JT85}]
\label{prop_JT}
Let $({\cal M}, \hat{g})$ be an ASD four-manifold with a conformal Killing
vector $K$.
The EW structure on the space  ${\cal W}$ of trajectories of $K$ 
{\em (}which is assumed to be non-pathological{\em)} is
defined by
\be
\label{EWs}
h:=|K|^{-2}\hat{g}-|K|^{-4}{\bf K}\odot {\bf K},\;\;\; \nu:=s^*\left(2|K|^{-2}
\ast_{\hat{g}}({\bf K}\wedge \d{\bf K})\right),
\ee
where $|K|^2:=\hat{g}_{ab}K^aK^{b}$, ${\bf K}$ is the one form dual to $K$ and
$\ast_{\hat{g}}$ is taken with respect to $\hat{g}$ and $s:{\cal
W}\mapsto {\cal M}$ is an arbitrary section of the fibration ${\cal M}\mapsto
{\cal W}$.  
All EW structures arise in this way.

Conversely, let $(h, \nu)$ be a three--dimensional EW structure 
on ${\cal W}$, and
let $(V, \a)$ be a pair consisting of a function of weight $-1$ 
and a one-form on ${\cal W}$ which satisfy
the generalised monopole equation
\be
\label{EWmonopole}
\ast_h(\d V+(1/2)\nu V) =\d\a,
\ee
where $\ast_h$ is taken with respect to $h$. Then
\be
\label{Vag}
g=Vh\pm V^{-1}(\d z+\a)^2
\ee
is an ASD metric with an isometry $K=\p_z$.
The minus sign in {\em{(\ref{Vag})}} is choosen if $h$ has signature
$(++-)$.
\end{prop}

In what follows we shall consider ASD structures which are 
also (complexified) hyper-Hermitian.

A smooth manifold ${\cal M}$ equipped with three almost complex
structures $(I, J, K)$ satisfying the algebra of quaternions is called
hyper-complex iff the almost complex structure $
{\cal J}_{\l}=aI+bJ+cK
$
is integrable for any $(a, b, c)\in S^2$.  
We use 
$\l=(a+ib)/(c-1)$, a stereographic coordinate on $S^2$ which we view
as a complex projective line $\CP^1$.
Let $g$ be a
Riemannian metric on $\cal M$.  If $({\cal M}, {\cal J}_{\l})$ is
hyper-complex and $g({\cal J}_{\l}X, {\cal J}_{\l}Y)=g(X, Y)$ for all
vectors $X, Y$ on ${\cal M}$ then the triple $({\cal M},\; {\cal
  J}_{\l},g)$ is called a 
hyper-Hermitian structure.  

We will in practice 
be interested in complexified or indefinite
hyper-Hermitian metrics with signature $(++--)$ for which the tensors
$(I,J,K)$ must necessarily be complex.  In signature $(++--)$ we can
arrange that one be real and the other two be pure imaginary, in the
latter case they determine a pair of transverse null foliations. 

We
shall restrict ourselves to oriented four manifolds.  In four
dimensions a hyper-complex structure defines a conformal structure,
which in explicit terms is represented by a conformal orthonormal
frame of vector 
fields $(X, IX, JX, KX)$, for any $X\in T{\cal M}$.
It is well known \cite{B88} that this conformal structure is ASD
with the orientation determined by the complex structures.

If there exists a choice of a conformal factor such that a two form
$\Sm_\l$ defined by $\Sm_\l(X, Y):=g(X, {\cal J}_{\l}Y)$ is closed (with
fixed $\l$) 
for  all $\l\in \CP^1$ and all vectors $(X, Y)$ then  $({\cal M},\; {\cal
J}_{\l},g)$ is called hyper-K{\"a}hler.

We shall use the following characterisation of the hyper-Hermiticity condition:
\begin{prop}[\cite{MN89,D98}]
\label{Laxmacia}
Let $\nabla_{AA'}$ be four independent real vector fields on 
a four-dimensional real manifold ${\cal M}$, and
let 
\[
L_0=\nabla_{00'}-\lambda\nabla_{01'},\;\;\;\;\;
L_1=\nabla_{10'}-\lambda\nabla_{11'},\;\;\;{\mbox{where}}\;\; 
\lambda \in \CP^1.
\]
If 
\be
\label{Lemat}
[L_0, L_1]=0
\ee
for every $\lambda$, then $\nabla_{AA'}$ is a null tetrad for a $(++--)$
hyper-Hermitian metric on ${\cal M}$. Every $(++--)$ hyper-Hermitian metric
arises in this way. Moreover, if the vectors $\nabla_{AA'}$ preserve
a volume form ${\mbox{{\em vol}}}_g$ on ${\cal M}$, then  $f^{-1}\nabla_{AA'}$ is a null tetrad
for a $(++--)$ hyper-K\"ahler metric on ${\cal M}$. Here
$f^2=\mbox{{\em vol}}_g(\nabla_{00'}, \nabla_{10'}, \nabla_{01'},
\nabla_{11'})$.
\end{prop}
Now we shall use (\ref{EWdkp}) and Proposition \ref{prop_JT} to construct  
ASD metrics out of solutions
to the dKP equation, and Proposition \ref{Laxmacia} to show that they are 
hyper-K{\"a}hler.

Assume that $h$ and $\nu$ are as in (\ref{EWdkp}).
Taking the exterior derivative of the generalised monopole equation
(\ref{EWmonopole}) yields
\begin{eqnarray}
\label{lindKP}
0&=&\nabla_i\nabla^iV+(1/2)(\nabla^i\nu_i)V+(1/2)\nu^i\nabla_iV\nonumber\\
&=&V_{yy}-V_{xt}+uV_{xx}+2u_xV_x+u_{xx}V
\end{eqnarray}
which is just a linearisation of the dKP equation (\ref{dKP})
(note that for $u=0$ (\ref{lindKP}) is just the wave equation relative
to the flat metric $\d y^2-4\d x\d t$).
One solution is $V=u_x/2$. One could find a corresponding $\a$ and write down
a metric using formula (\ref{Vag}) (see the remarks after 
Proposition \ref{prop_KCdKP}),
but we shall present a different method based on
the Lax operators.

Take the Lax operators (\ref{dKPlax}) and introduce a new spectral parameter
$\l:=\lt-z$ for some $z$. The function $u(x,y,t)$ does not depend on $z$ so we can
replace $\p_{\lt}$ by $\p_z$. This yields (with dropped primes and added tildes)
\begin{eqnarray*}
\tilde{L}_0&=&\p_t-u\p_x-z\p_y+u_y\p_z-\l \p_y,\\
\tilde{L}_1&=&\p_y-z\p_x+u_x\p_z-\l \p_x.
\end{eqnarray*}
To obtain a pair of exactly commuting operators take
\begin{eqnarray*}
{L}_1&:=&\t{L}_1=\p_y-z\p_x+u_x\p_z-\l \p_x,\\
{L}_0&:=&\t{L}_0+z\t{L}_1=\p_t-(u+z^2)\p_x+(u_y+u_xz)\p_z-\l(\p_y+z\p_x). 
\end{eqnarray*}
If $u(x,y,t)$ is a solution to (\ref{dKP}) then these operators 
satisfy $[{L}_0, L_1]=0$ and so, by Proposition \ref{Laxmacia}, the vectors 
\[
\nabla_{10'}=\p_y-z\p_x+u_x\p_z,\;\nabla_{11'}=\p_x,\;
\nabla_{00'}=\p_t-(u+z^2)\p_x+(u_y+u_xz)\p_z,\;\nabla_{01'}=(\p_y+z\p_x),
\]
form a hyper-Hermitian frame. The vectors $\nabla_{AA'}$ preserve the
volume form  $\mbox{vol}_{g}=\d t\wedge\d
y\wedge \d x \wedge \d z$, and $f^2=u_x/2$. 
Therefore we have the following
\begin{prop}
\label{prop_KCdKP}
Let $u=u(x, y, t)$. The metric
\be
\label{HCdKP}
g=\frac{u_x}{2}(\d y^2-4\d x\d t-4u\d t^2)
-\frac{2}{u_x}(\d z-\frac{u_x\d y}{2}-u_y\d t)^2
\ee
is hyper-K{\"a}hler.
\end{prop}
{\bf Remarks:}
\begin{itemize}
\item The above metric has a Killing vector $\p_z$ with the dual
\[
K=-\frac{2}{u_x}(\d z-\frac{u_x \d y}{2}-u_y\d t),
\]
and the formulae (\ref{EWs}) gives rise to the Einstein--Weyl structure
(\ref{EWdkp}). 
The self-dual part of $\d K$ is a simple two-form. In section
\ref{twistor_theory} we shall show that all hyper-K{\"a}hler metrics
with such symmetries are locally given by (\ref{HCdKP}).
\item Note that $u_x\neq 0$ for (\ref{HCdKP}) to be well defined. 
To obtain a flat metric take $u=-x/t$
which is a special case of (\ref{conformal_Einstein}). 
The metric (\ref{HCdKP})
becomes
\[
g=2\d x\frac{\d t}{t}-2x\frac{\d t^2}{t^2}+2t\d z^2+2\d z\d y.
\]
Putting $
x=Xt+{z^2t}/2,\;y=Y-zt $ yields the flat metric
\[
g=2\d X\d t+2\d z\d Y.
\]

\item The metric (\ref{HCdKP}) could be found directly from the
monopole equation (\ref{EWmonopole}) as follows: 
Rewrite the metric (\ref{EWdkp}) in an orthonormal triad
$h=e_1^2+e_2^2-e_3^2$, where
\[
e_1=\d y,\qquad e_2=\d x+(u-1)\d t,\qquad e_3=\d x+(u+1)\d t.
\]
The duality relations $
\ast_h e_1=e_3\wedge e_2,\; \ast_h e_2=e_1\wedge e_3
,\; \ast_h e_3=e_1\wedge e_2$
yield
\be
\label{du_rel}
\ast_h\d t=\d t\wedge\d y,\qquad \ast_h\d y=2\d t\wedge\d x,\qquad 
\ast_h\d x=\d y\wedge\d x+2u\d y\wedge\d t.
\ee
Take $V=u_x/2$, and use the above relations to
write the monopole equation (\ref{EWmonopole}) as
\[
\frac{u_{xx}}{2}\d y\wedge\d x +u_{xy}\d t\wedge\d x
+(u_x^2+uu_{xx}-\frac{u_{xt}}{2})\d y\wedge\d t
=\d \a.
\]
Choosing the gauge in which $\a= \a_1\d y+\a_2\d t$ 
(this is always possible by redefining a coordinate $z$ along the
orbits of a Killing vector) gives 
\be
\label{eqn_for_a}
(\a_1)_x=-\frac{u_{xx}}{2},\qquad (\a_2)_x=-u_{xy},\qquad
(\a_2)_y-(\a_1)_t=\frac{u_{xt}}{2}-u_{yy}.
\ee
All solutions to this system of equations are gauge equivalent to
\[
\a=-\frac{u_x}{2}\d y-u_y\d t.
\]
Substituting $V, \a$ and $h$ to (\ref{Vag}) yields (\ref{HCdKP}).
\item The Lax pair (\ref{dKPlax}) can be obtained from the 
hyper-K{\"a}hler Lax pair by a symmetry reduction:
The distribution
$(K, \t{L}_0, \t{L}_1)$ is not integrable, as $[K, \t{L}_0]=-\p_y$ and
$[K, \t{L}_1]=-\p_x$. To obtain an integrable distribution, one needs to lift $K$ to the correspondence space by
$\tilde{K}=K-\p_{\l}$. Then $(\tilde{K}, \t{L}_0, \t{L}_1)$ is an integrable distribution,
but $\tilde{K}(\l)\neq 0$, which forces us to introduce an  invariant spectral parameter
$\lt=\l+z$. This implies that in the Lax pair we replace all $\p_z$ by $\tilde{K}+\p_{\lt}$. Now we restrict ourselves to invariant solutions to 
$\t{L}_0\Psi=\t{L}_1\Psi=0$,
and so we ignore $\Kt$ in the Lax pair. The reduced Lax pair is  given by
(\ref{dKPlax}).
\end{itemize}
In the covariantly constant primed spin frame the null tetrad is
\begin{eqnarray*}
e^{00'}&=&-u_x\d t,\qquad e^{10'}=\frac{\d z-u_y\d t}{u_x},\\
e^{01'}&=&\d z-u_x\d y-(u_y+zu_x)\d t,\qquad e^{11'}=
\d x+u\d t+z\frac{\d z-u_y\d t}{u_x},
\end{eqnarray*}
and the metric (\ref{HCdKP})  is $2(e^{00'}e^{11'}-e^{01'}e^{10'})$.
The basis of SD two form is in this frame given by
\begin{eqnarray*}
\Sm^{0'0'}&=&\d z\wedge\d t,\qquad
\Sm^{0'1'}=\d z\wedge\d y+\d (u+z^2)\wedge\d t,\\
\Sm^{1'1'}&=&u_x\d x \wedge\d y-uu_x\d y\wedge \d t+u_y\d x\wedge \d t
+\d(uz)\wedge\d t+\d z\wedge(\d x+z\d y+z^2\d t).
\end{eqnarray*}
They satisfy 
\[
-2\Sm^{0'0'}\wedge\Sm^{1'1'}=\Sm^{0'1'}\wedge\Sm^{0'1'},
\qquad \d \Sm^{0'0'}=\d \Sm^{0'1'}=\d \Sm^{1'1'}=0,
\]
which again implies that the metric (\ref{HCdKP}) is hyper--K\"ahler. 
Note that the Killing vector $K=\p_z$ does not preserve the K{\"a}hler
form $\Sm^{0'1'}$.




\section{Examples}
\subsection{dKP EW spaces with $S^1$ symmetry}
In this subsection we shall construct  EW structures depending on one 
arbitrary function
of one variable.

To find some explicit examples of (\ref{EWdkp})
assume that $u$ is independent of $y$. Therefore it satisfies the 
simple equation
$uu_x=u_t$, all solutions of which are given in an implicit form 
\[
u(x,t)=f(x+tu(x,t))
\]
(more general hodograph transformations for dKP 
arising from its connection with equations of hydrodynamic type
were studied in \cite{K88}, and \cite{GK89}).  
 
Here $f$ is an arbitrary function of one variable $s:=x+tu(x,t)$. The idea is to 
write the Einstein--Weyl structure (\ref{EWdkp}) 
making use of this `hodograph transformation'.
We have
\[
h=\d y^2-4\d t(\d x+ u\d t)=\d y^2-4\d t(\d s-t\d u)=\d y^2 -4\d t\d s+4t\d t\d f(s)
\]
where we performed a coordinate 
transformation $(x, y, t)\rightarrow (s, y, t)$.
Defining $F(s):=\d f/\d s$ and replacing $u_x$ by $F/(1-tF)$ yields
the EW structure
\be
\label{dKdV}
h=\d y^2+4(tF(s)-1)\d t\d s,\;\;\;\;\;\;\; \nu=4\frac{F(s)}{tF(s)-1}\d t,
\ee
which depends on one arbitrary function $F(s)$
 (which we shall take to be strictly negative) of one variable. 
This structure has  signature $(++-)$. If $t>0$ then it is
well-defined on $S^1\times \R^+\times \R$.

We shall now show that
formulae (\ref{dKdV}) give a class of $EW$ structures  on principal
$S^1$ bundles over Weyl manifolds.
\begin{prop}
\label{Mexico}
Let $({\cal N}, [H], \nu_H)$ be a two-dimensional manifold 
with a Weyl structure of signature $(+-)$ 
and let $\pi:{\cal W}\longrightarrow {\cal N}$
be an $S^1$ bundle over $N$.  
If
\[
h:=\d y^2 +\pi^*H,\;\;\;\;\; \nu:=\pi^*\nu_H
\]
(where $y$ is a coordinate on a fibre) is an EW structure on ${\cal
W}$ then it can be put in  the form \em{(\ref{dKdV})}. 
\end{prop}
{\bf Proof.}
We can use isothermal coordinates $(\tilde{s}, t)$ on ${\cal N}$ and choose a
representative of a conformal class $[H]$ such that  $h$ and $\nu$ are
\be
\label{Qgen}
h=\d y^2+2G(\tilde{s},t)\d \tilde{s}\d t,\;\;\;\;\;\;\; \nu=K(\tilde{s},t)\d t.
\ee
Each EW structure of this form is equivalent to (\ref{dKdV}). This can be seen as follows: Equations
$\chi_{13}=0, \chi_{22}=0$ imply that $K=4G_t/G+f(t)$. The function $f(t)$ can be
absorbed in the definition of $G$. Then  the vanishing of $\chi_{33}$
(all remaining EW equations are satisfied trivially) yields 
$G(\tilde{s},t)=-2F_1(\tilde{s})+2tF_2(\tilde{s})$ 
for arbitrary $F_1$ and $F_2$. Now we define a new coordinate $s$ by 
$\d s:= F_1(\tilde{s})\d \tilde{s}$. Equivalence between (\ref{Qgen}) and
(\ref{dKdV}) is finally obtained by putting $F(s):=F_2(s)/F_1(s)$.
The metric (\ref{dKdV}) is not Einstein as $G_{22}\neq 0, G_{13}\neq 0$
and $R=-2F_s/(tF-1)^3$ is not constant (unless $F$ is constant).
To visualise the two-dimensional surface $\cal N$ on which $H$ is defined  
one can restrict a flat $(++--)$ metric on $\R^{4}$, $g= \d f\d w-\d s\d t$ to the 
intersection of the paraboloid $w=t^2/2$ with the hyper-surface $f=f(s)$.
\koniec

The hyper-K{\"a}hler metric corresponding to (\ref{dKdV}) 
has an additional null Killing vector $\p_y$ and is  
(with definitions $\d w:=-F\d s, \hat{F}(w):=F^{-1}$) given by
\[
g=\d w\d t+\d z\d y +(t-\hat{F}(w))\d z^2
\]
where $\hat{F}(w)$ is arbitrary.

 Other examples (without a Killing vector) can be obtained from
\[
u=t\frac{\d A(t)}{\d t}-\frac{x}{t}+\frac{y}{t}\sqrt{\frac{x}{t}+A(t)},
\]
where $A(t)$ is arbitrary.

\subsection{dKP metrics which are hyper-CR}  
Let us recall that that an EW metric is called hyper-CR (or special) 
if it admits a two-parameter family of
shear-free, divergence-free geodesic congruences \cite{CT00}. All hyper-CR EW
spaces arise as reductions of hyper-K\"ahler metrics by triholomorphic 
homotheties \cite{GT98}. 
In this section we shall find all EW metrics in 2+1 
dimensions which are both dKP and hyper--CR.
This will lead to a class of solutions to the dKP equation
depending on one arbitrary function of one variable.

\begin{prop}
All EW metrics which admit a constant weighted vector and a 
two parameter family of shear-free geodesic 
congruences with a vanishing divergence 
are either spaces of constant curvature or are locally of the
form
\be
\label{dKP_special}
h=\d y^2-4\d x\d t-4\Big(\frac{P(t)}{y}-\frac{x^2}{y^2}\Big)\d t^2,
\qquad \nu=\frac{8x}{y^2}\d t,
\ee
where $P$ is an arbitrary function of $t$.
\end{prop}
{\bf Proof.} 
The hyper-CR condition for a metric is characterised \cite{GT98} 
by the existence of a scalar $\rho$ of weight $-1$ which (together
with the Einstein--Weyl one form $\nu$)  
satisfies the monopole equation 
\be
\ast_h(\d \rho +\frac{1}{2}\nu\rho)=\d \nu,
\ee
and the algebraic constraint
\be
\label{special}
\qquad \rho^2=\frac{8}{3}W.
\ee
We shall impose these conditions on the dKP metric (\ref{EWdkp}).
The monopole  equation yields
\[
(4u_{xx}-2\rho_y)\d x\wedge\d t+\rho_x\d y\wedge\d x+
(2\rho_x u-\rho_t+2\rho u_x+4u_{xy})\d y\wedge\d t=0
\]
which (together with (\ref{special})) gives four scalar equations:
\be
\label{set_of _eq}
\rho_y=2u_{xx}, \qquad \rho_x=0, \qquad 2\rho u_x-\rho_t+4u_{xy}=0,
\qquad \rho^2=-8u_{xx}.
\ee
If $u_{xx}=0$ then the last relation in 
(\ref{set_of _eq}) gives  $\rho=0$. The monopole equation  then implies
that  $\nu$ is closed, and the Einstein--Weyl metric is conformal to
Einstein. Therefore we assume $u_{xx}\neq 0$. Differentiating the 
third equation in (\ref{set_of _eq}) with respect to $x$ (and
using the first two equations) gives 
\[
\rho=-2\frac{u_{xxy}}{u_{xx}}.  
\]
The integrability
conditions to (the otherwise over-determined system) (\ref{set_of _eq})
are
\begin{eqnarray}
\label{integrability_con}
& &u_{xxx}=0,\qquad u_{xxy}^2-u_{xxyy}u_{xx}=u_{xx}^3,\qquad
4u_{xxy}=\eta u_{xx}^3,\\
& &u_{xxy}u_{xxt}-u_{xxyt}u_{xx}+
2u_xu_{xx}u_{xxy}-2u_{xy}u_{xx}^2=0.\nonumber
\end{eqnarray}
The first condition implies
$u(x,y,t)=ax^2+bx+c$. Here $a, b, c$ are functions of $y$ and $t$,
which satisfy
\be
\label{int_1}
a_{yy}+6a^2=0,
\ee
\be
\label{int_2}
b_{yy}-2a_t+6ab=0,
\ee
\be
\label{int_3}
c_{yy}-b_t+2ac+b^2=0,
\ee
\be
\label{int_4}
a_{y}^2-aa_{yy}-2a^3=0,
\ee
\be
\label{int_5}
a_{y}^2+4a^3=0,
\ee
\be
\label{int_6}
aa_{yt}-a_ya_t-2aa_yb+2b_ya^2=0
\ee
Equations (\ref{int_1}, \ref{int_2}, \ref{int_3}) follow from 
the dKP (\ref{dKP}), and the other equations are the integrability 
conditions (\ref{integrability_con}). Solve  (\ref{int_5}) to find
$
a(y, t)=-(y-L(t))^{-2}
$
(or $a=0$ which gives $u_{xx}=0$).

We can now perform the coordinate transformation (\ref{coordtrans})
with $f=-L/2$ and $g=0$
to set $L(t)=0$.
One verifies that  (\ref{int_1}), and (\ref{int_5}) are now also satisfied.
Equation (\ref{int_2}) gives
$
b(y, t)=-M(t)y^{-2}+N(t)y^3,
$
but (\ref{int_6}) implies $N(t)=0$. 
So far we have
\[
h=\d y^2-4\d x \d t+4\Big(c(y, t)-\frac{xM(t)}{y^2}-\frac{x^2}{y^2}\Big)\d t^2,
\qquad
\nu=\frac{8x+4M(t)}{y^2}\d t.
\]
The function $M(t)$ can  be eliminated by the coordinate
transformation (\ref{coordtrans}) with $g=M/2$. 
Imposing (\ref{int_3}) yields  
$c(y, t)=P(t)/y + R(t)y^2$ leaving
\[
h=\d y^2-4\d x \d
t+4\Big(-\frac{x^2}{y^2}+\frac{P(t)}{y}+R(t)y^2\Big)\d t^2,\qquad
\nu=\frac{8x}{y^2}\d t.
\]


We eliminate $R(t)$ by performing the conformal rescaling and associated
coordinate redefinitions of (\ref{conftrans}) with $c(\tilde{t})$ satisfying
$$
R=-\frac{c^{\prime\prime\prime}}{6c^{\prime 3}}
+\frac{1}{4}\left(\frac{c''}{c^{\prime 2}}\right)^2\, .
$$
This yields, dropping the tildes and with a redefinition of $P$,
\[
u(x, y, t)=-\frac{x^2}{y^2}+\frac{P(t)}{y}.   
\]
The Einstein--Weyl structure is therefore 
(\ref{dKP_special}).
The arbitrary function $P(t)$ can not be
eliminated. 
This can be seen by finding the symmetries (\ref{Weyl_Killing})
of the EW structure (\ref{dKP_special}). We summarise our findings in
the table below:
{\small
\begin{center}
\begin{tabular}{p{1cm}|lll}
\multicolumn{3}{c}{}\\ 
 &Function $P(t)$&Symmetries\\
\hline\
{\bf (i)}&$P(t)=0$& $ K_1, K_2, K_3, K_4$\\
{\bf (ii)}&$P(t)=const\neq 0$&  $K_1, K_2+3K_3, K_4$\\
{\bf (iii)}&$P(t)=(bt+c)^{\frac{3a-b}{2b}}$& $cK_1+aK_2+bK_3$\\
{\bf (iv)}&general $P(t)$&none
\end{tabular}
\end{center}
\small}
\noindent
where $a, b, c$ are constants, and 
\[
K_1=\p_t,\qquad K_2=(1/2)y\p_y+x\p_x,\qquad K_3=(1/2)y\p_y+t\p_t,\qquad
K_4=ty\p_y+(y^2+2xt)\p_x+3t^2\p_t.
\]
Note that in case $(ii)$ we can redefine coordinates to set $P(t)=1$.
The vector fields 
$K_1, K_2+3K_3, K_4$ generate the Lie group of Bianchi type VIII,
i.e.\ $SU(1,1)$,
and the cases $(i)$ and $(ii)$ give homogeneous EW spaces.
Case $(iii)$ can be reduced to $P(t)=t^{\a}, K=K_3+[(2\a+1)/3]K_2$,
where $\a=const\neq 0$.
\koniec

\section{The twistor correspondences and Lax formulations}
In this section we shall study the twistor theory of the EW spaces.
We first discuss the twistor
correspondence in the flat case.  We then give a Lax
formulation of the EW equations and derive from it the twistor
correspondence. 
We study this correspondence in relation to reductions of the
anti-self-duality equations on four-dimensional conformal structures.   
We then reformulate the Einstein--Weyl equations in terms
of a certain two-form on the trivial $\CP^1$ bundle over a Weyl space.

\subsection{The flat correspondence}
\label{EW_bundle}  
Let us begin by recalling Ward's approach \cite{W89} to twistors in 
(2+1)-dimensional
flat space-times.
Rearrange the space time coordinates $(x, y, t)$ as a symmetric 
two-spinor\footnote{The use of primed (rather than unprimed) spinors
in this section originates from the representation of
Einstein--Weyl spaces as reductions of ASD (rather than SD) metrics in
four dimensions.
ASD structures (for which the covariantly constant self-dual spinors
are conventionally denoted as having primed indices)
are taken as basic because they arise from a natural
choice of orientation and conformal structure on a K\"ahler manifold. 
}
\[
x^{A'B'}:=
\left (
\begin{array}{cc}
t&y/2\\
y/2&x
\end{array}
\right ),
\]
such that the space-time metric and the volume form are :
\[
h=-2\d x_{A'B'}\d x^{A'B'},\qquad 
\mbox{vol}_h=\d {x_{A'}}^{B'}\wedge\d {x_{C'}}^{A'}\wedge\d {x_{B'}}^{C'}.
\]
The two-dimensional spinor indices are raised and lowered with the symplectic form 
$\varepsilon_{A'B'}$, such that $\varepsilon_{0'1'}=1$
(see \cite{PR86} for a full account of the two-spinor formalism).
We shall use the abstract index convention
$V^i=V^{(A'B')}=v^{(A'}\pi^{B')}$ based on an isomorphism
$T^i{\cal W}=S^{(A'}\otimes S^{B')}$.

The projective mini-twistor space of $\R^{2+1}$ is the two-dimensional complex manifold 
${\cal Z}=T\CP^1$ which is the total space of the line bundle 
${\cal O}(2)$ 
of Chern class 2 over $\CP^1$.
Points of ${\cal Z}$ correspond to null 2-planes in $\R^{2+1}$ via the 
incidence relation
\be
\label{minit1}
x^{A'B'}\pi_{A'}\pi_{B'}=\om.
\ee
Here $(\om, \pi_{0'}, \pi_{1'})$ are homogeneous coordinates on ${\cal
O}(2)$: $(\om, \pi_{A'})\sim (\rho^2 \om, \rho\pi_{A'})$, 
where $\rho \in \C^*$.
 In the affine coordinates $\lt:=\pi_{0'}/\pi_{1'}, 
\xi:=\om/(\pi_{1'})^2$
equation (\ref{minit1}) is
$
\xi=x+\lt y+\lt^2 t.
$
First fix $(\om , \pi_{A'})$. If $(\xi, \lt)$ are both
real then (\ref{minit1})  defines a null plane in $\R^{2+1}$. If both $\xi$ and $\lt$ are complex then
the solution to (\ref{minit1})  is a time like curve in $\R^{2+1}$. We shall say that this curve is
oriented to the future if $\mbox{Im}{\lt}>0$ and to the past
otherwise. 
If $\lt$ is real and $\xi$ is complex then
(\ref{minit1}) has no solutions for finite $x^{A'B'}$.

An alternate interpretation of (\ref{minit1}) is to fix $x^{A'B'}$.
This determines $\om$ as a function of $\pi_{A'}$ i.e.\ a section of
${\cal O}(2)\rightarrow \CP^1$ when factored out by the relation
$(\om, \pi_{A'})\sim(\rho^2 \om, \rho\pi_{A'})$. 
These are embedded rational curves with normal bundle 
$\O(2)$. Two rational curves $l_{p_1}$ and $l_{p_2}$ (corresponding
to $(t_1, y_1, x_1)$ and $(t_2, y_2, x_2)$ respectively) intersect at two points
\[
\l_{1, 2}=\frac{2R_2\mp\sqrt{h(R, R)}}{2R_1},\qquad 
\mbox{where}\qquad R_i:=(t_1-t_2, y_1-y_2, x_1-x_2).
\]
Therefore the incidence of curves in ${\cal Z}$ encodes the causal
structure of $\R^{2+1}$ in the following sense: $l_{p_1}$ and
$l_{p_2}$ intersect at {\bf (a)} one point, {\bf (b)} two real points,
{\bf (c)} two complex points conjugates of each other, iff $p_1, p_2$
are {\bf (a)} null separated, {\bf (b)} space-like separated, {\bf
  (c)} time-like separated.

Examining the relevant cohomology groups shows that the moduli space
of curves with normal bundle $\O(2)$ in $\cal Z$ is $\C^3$. The real
space-time $\R^{2+1}$ arises as the moduli space of curves that are
invariant under the conjugation $(\omega,\pi_{A'})\mapsto
(\bar{\omega},\bar{\pi}_{A'})$.

The correspondence space ${\cal F}=\C^3\times\CP^1
=\{(p,Z)\in\C^3\times {\cal Z}|{Z\in l_p}\}$.  By definition, it
inherits fibrations over both $\C^3$ and ${\cal Z}$ and the fibration
of ${\cal F}=\C^3\times\CP^1$ over ${\cal Z}$ has fibres spanned by
the distribution $ L_{A'}=\pi^{B'}\p_{A'B'}$, where $\p_{A'B'}x^{C'D'}
=1/2(\varepsilon^{C'}_{A'}\varepsilon^{D'}_{B'}
+\varepsilon^{C'}_{B'}\varepsilon^{D'}_{A'})$.  In the affine
coordinates $\pi^{A'}=(1, -\lt)$ this distribution is
\[
L_{0'}=\p_t-\lt\p_y,\;  L_{1'}=\p_y-\lt\p_x
\]
(we have ignored the constant factor $\pi_{1'}$). Note that this
$L_{A'}$ is the special case $u(x, y, t)=0$ of the Lax pair
(\ref{dKPlax}) for the dKP equation.

We also define the correspondence space ${\cal F}_W=\R^{2+1}\times
\CP^1$ for $\R^{2+1}$. 
Let ${\cal Z}_{\R}$ be the sub-manifold of ${\cal Z}$ preserved by the
conjugation
\[
(\om, \pi_{0'}, \pi_{1'})\rightarrow(\ov{\om}, \ov{\pi_{0'}},
\ov{\pi_{1'}}),
\]
and let $l_p$ be the real line in ${\cal Z}_{\R}$ 
that corresponds to $p\in {\cal W}$ and let $Z\in l_p$.
The totally real correspondence space is a four-dimensional real manifold defined by
${\cal F}^4_\R:={\cal Z}_\R\times \R^{2+1}|_{Z\in l_p}$ and can be
represented as the set $\lt=\bar{\lt}$ or $\pi_{A'}=\bar{\pi}_{A'}$. The distribution
$L_{A'}\cap \overline{L}_{A'}$ is one dimensional, spanned by
$\ov{\pi}^{A'}\pi^{B'}\p_{A'B'}$, on the complement of ${\cal F}_\R^4$. On ${\cal F}^4_\R$ 
$L_{A'}\cap \overline{L}_{A'}$ is two real dimensional as here 
$L_{A'}=\bar{L}_{A'}$. 
The real correspondence space ${\cal F}_{\R}$ divides ${\cal
  F}_W=\R^{2+1}\times \CP^1$ into two halves.

\subsection{The Lax formulation and twistor correspondence}
\begin{prop}
\label{LAX_PAIR_THEOREM}
Let $V_1, V_2, V_3$ be three independent holomorphic
vector fields on a three dimensional complex manifold ${\cal W}$ 
such that
\be
\label{LAX_PAIR}
L_{0'}=V_1-\lt V_2+f_{0'}\p_{\lt},\qquad
L_{1'}=V_2-\lt V_3+f_{1'}\p_{\lt}
\ee
is an integrable distribution for some functions $f_{0'}, f_{1'}$, 
which are third-order polynomials in $\lt\in\CP^1$. Then there exists a one form 
$\nu$  such that the contravariant metric
$V_2\otimes V_2-1/2(V_1\otimes V_3+ V_3\otimes V_1)$ 
and $\nu$ give an EW structure on ${\cal W}$. 
Each EW structure
arises in this way.
\end{prop}
{\bf Remarks:}
\begin{itemize}
\item The Lax pair (\ref{dKPlax}) for the dKP equation is of course a 
special case of (\ref{LAX_PAIR}).
\item The  Lax formulations are widely applicable in the theory of integrable
systems and so the above proposition
can be applied outside  twistor theory. 
It is however much easier to prove Proposition 
\ref{LAX_PAIR_THEOREM} using the twistor geometry, rather
than an explicit calculation. This justifies adopting the spinor notation
\[
\nabla_{A'B'}=
\left (
\begin{array}{cc}
V_1&V_2\\
V_2&V_3
\end{array}
\right ),
\qquad f_{A'}=(f_{0'}, f_{1'}), \qquad \pi^{A'}=(1, -\lt),
\]
in which the Lax pair has the compact form
$L_{A'}=\pi^{B'}\nabla_{A'B'}+f_{A'}\p_{\lt}$.
We shall use this notation in the proof of Proposition 
\ref{LAX_PAIR_THEOREM}.
\item The third order polynomials $f_{A'}$ contain eight functions 
not depending on $\lt$. These can be reduced to four functions 
by choice of a suitable spin frame for which
$f_{A'}$  become linear in $\lt$.
In this  frame
there exists a vector formula for $\nu$ in terms of $\Gamma_{ijk}$,
and $f_{A'}$.
\item  Proposition \ref{LAX_PAIR_THEOREM} holds for complex solutions
  and for any choice of signature for real space time.
\end{itemize}
{\bf Proof of Proposition \ref{LAX_PAIR_THEOREM}.} 
Assume that $h=V_2\otimes V_2-1/2(V_1\otimes V_3+ V_3\otimes V_1)$ 
and $\nu$ gives an EW structure.
Let $V(\lt)=V_1-2\lt V_2+\lt^2V_3$. Then $g(V(\lt), V(\lt))=0$ for 
all $\lt\in\CP^1$
so $V(\lt)$ determines a sphere of null vectors. 
Choose  
$l_{0'}=V_1-\lt V_2,\;
l_{1'}=V_2-\lt V_3$ as a basis of the orthogonal complement
of $V(\lt)$. For each $\lt\in\CP^1$ the vectors $l_{0'}, l_{1'}$ give
a null two-surface. It is well known \cite{C43, H82, PT93} 
that the EW equations on $(h, \nu)$ are equivalent to 
the integrability conditions of null, totally geodesic planes.
Therefore  the Frobenius theorem implies that the 
horizontal lifts 
\[
L_{0'}=V_1-\lt V_2+f_{0'}\p_{\lt},\qquad
L_{1'}=V_2-\lt V_3+f_{1'}\p_{\lt}
\]
of $l_{0'}, l_{1'}$ to $T({\cal W}\times\CP^1)$
span an integrable  distribution. The functions $f_{0'}$ and  $f_{1'}$ 
are third order in $\lt$, because the M{\"o}bius transformations of
$\CP^1$ are generated by  vector fields quadratic in $\lt$, 
and  $l_{0'}, l_{1'}$ are linear $\lt$.

The above argument can be made more explicit in spinor notation:
let $L_{A'}$ be horizontal lift of $l_{A'}=\pi^{B'}\nabla_{A'B'}$
to the weighted spin bundle (i.e. $L_{A'}\pi_{C'}=0$).
This yields
\begin{eqnarray}
\label{EWlax}
L_{A'}&=&\pi^{B'}\nabla_{A'B'}+
\Gamma_{A'B'C'D'}\pi^{B'}\pi^{D'}\frac{\p}{\p \pi_{C'}}\nonumber\\
& &+\frac{1}{2}\nu_{B'D'}\pi^{B'}\Big(\pi^{D'}\frac{\p}{\p \pi^{A'}}
-\frac{1}{2}\pi_{A'}\frac{\p}{\p \pi_{D'}}-
{\varepsilon_{A'}}^{D'}\pi\cdot\frac{\p}{\p \pi}\Big),
\end{eqnarray}
where $\Gamma_{A'B'C'D'}$ is spinor Levi--Civita connection defined
by $\nabla_{A'B'}\pi_{C'}=-\Gamma_{A'B'C'D'}\pi^{D'}$.
The integrability conditions imply
$
[L_{A'}, L_{B'}]= 0\;(\mbox{mod}\;L_{A'}).
$
The distribution $L_{A'}$, when projected to $\cal F_W$ is given by
(\ref{LAX_PAIR}), where 
\[
f_{A'}=\Gamma_{A'B'C'D'}\pi^{B'}\pi^{C'}\pi^{D'}+
(1/4)\pi_{A'}\nu_{B'C'}\pi^{B'}\pi^{C'}.
\]
\koniec

The twistor space ${\cal Z}$ for a solution to the  
EW equations on $({\cal W}, h, \nu)$ 
associated to the Lax system on $L_{A'}$ as above is obtained by
factoring the spin bundle ${\cal W}\times \CP^1$ by the twistor
distribution (Lax pair) ${L_{A'}}$. This clearly has a projection
$q :{\cal W} \times \CP^1 \mapsto {\cal Z}$ and we have a double
fibration
$$
\begin{array}{rcccl}
&&{\cal W}\times\CP^1&&\\
&r\swarrow&&\searrow q&\\
&{\cal W}&&{\cal Z}&
\end{array}
$$
Each point $p\in\cal W$ determines a sphere $l_p$ made up of all the
null totally geodesic two--surfaces through $p$. The normal bundle of $l_p$
in $\cal Z$ is $N=T{\cal Z}|_{l_p}/Tl_p$. This is a rank one vector
bundle over $\CP^1$, therefore it has to be one of the standard line
bundles ${\cal O}(n)$.
 \begin{lemma} The holomorphic curves $l_p:=q(\CP^1_p)$ where
$\CP^1_p=r^{-1}(p)$, $p\in {\cal W}$, have normal bundle $N={\cal
O}(2)$.  
\end{lemma}
\noindent
{\bf Proof.} To see this, note that $N$ can be identified with the
quotient $r^*(T_p{\cal W})/\{ \mathrm{span }\;L_{0'},L_{1'}\}$.
In their homogeneous form the operators $L_{A'}$ have weight one, so the
distribution spanned by them is isomorphic to the bundle
$\C^{2}\otimes{\cal O}(-1)$.  The definition of the normal bundle as
a quotient gives a sequence of sheaves over $\CP^1$.
\[
0\longrightarrow \C^{2}\otimes{\cal O}(-1) \longrightarrow \C^{3}  
\longrightarrow N\longrightarrow 0
\]
and we see that $N={\cal O}(2)$, because the
last map, in the spinor notation, 
is given explicitly by $V^{A'B'}\mapsto
V^{A'B'}\pi_{A'}\pi_{B'}$ clearly projecting
onto ${\cal O}(2)$.
\koniec
\smallskip

A generalisation of the flat mini-twistor correspondence 
to the 2+1  EW spaces is given by the following proposition
\begin{prop}[\cite{H82}]
Any solution to the EW equations {\em(\ref{ew2})} is equivalent
to a complex surface ${\cal Z}$  
with a family of rational curves with normal
bundle $\O(2)$.
\end{prop}
Points of ${\cal W}$ correspond to curves in ${\cal Z}$  with self-intersection
number 2. The Kodaira theorem \cite{Ko63} applied to deformations preserving the real 
structure of ${\cal Z}$
 guarantees the existence of a three-dimensional complex family
of such curves. Points of ${\cal Z}$ correspond to totally
geodesic hyper-surfaces in ${\cal W}$. Non-null geodesics in 
${\cal W}$ consist of all the curves in ${\cal Z}$ which
intersect at two fixed points in ${\cal Z}$. Null geodesics
correspond to curves passing through one point with a given tangent
direction.  Thus the projective and conformal structures can be
reconstructed. \koniec
\smallskip

\subsection{Mini-twistor spaces from twistor spaces}

\begin{prop}
All Einstein--Weyl spaces  arise as symmetry
reductions of hyper-Hermitian metrics (or indefinite 
hyper-Hermitian metrics)
in four-dimensions.
\end{prop} 
{\bf Proof.} Consider an EW structure with the corresponding Lax pair
(\ref{LAX_PAIR}).
Choose a spin frame in which $f_{A'}$ is linear in $\lt$;
$f_{A'}=U_{A'}+\lt W_{A'}$
(this is always possible by making a suitable M\"obius transformation
of $\CP^1$ and choosing an appropriate conformal scale) , and
introduce a new spectral parameter
$\l:=\lt-z$ for some $z$. Nothing in the $L_{A'}$ 
depends on $z$ so we can replace $\p_{\lt}$ by $\p_z$. This yields
(with a dropped prime)
\[ 
L_{A}=\nabla_{A0'}-\l\nabla_{A1'},
\]
where
\begin{eqnarray*}
\nabla_{00'}&=&\nabla_{0'0'}+z\nabla_{0'1'}+(U_{0'}+
 zW_{0'})\p_z,\\
\nabla_{10'}&=&\nabla_{1'0'}+z\nabla_{1'1'}+(U_{1'}+
 zW_{1'})\p_z,\\
\nabla_{01'}&=&\nabla_{0'1'} +W_{0'}\p_z, \\    
\nabla_{11'}&=&\nabla_{1'1'} +W_{1'}\p_z
\end{eqnarray*} 
where $U_{0'}, U_{1'}, W_{0'}, W_{1'} $ are four functions not
depending of $\l$. 
One is left with a Lax pair for a hyper-Hermitian four manifold
because  $L_{A}$ can be made to commute exactly
(as in  Proposition \ref{Laxmacia})
by choosing two solution to the background coupled
neutrino equation (see
\cite{D98} for details).
This Lax pair has an obvious symmetry $\p_z$. 
\koniec
\smallskip

\noindent 
{\bf Remark:}
All EW spaces arise as symmetry reductions
of a pair of coupled PDEs \cite{D98}, \cite{GS99} associated to 
hyper-Hermitian four manifolds.
In \cite{CP99} Proposition \ref{LAX_PAIR} was proven using different methods
for EW spaces of Riemannian signature.


The twistor construction of Hitchin can be viewed as a reduction of Penrose's
Nonlinear Graviton construction. It follows from \cite{JT85} (compare Proposition \ref{prop_JT})
that the mini-twistor space 
${\cal Z}$ corresponding to ${\cal W}$ is a factor space 
${\cal PT}/{\cal K}$ where
$\cal PT$ is the twistor space of $({\cal M}, g)$ and
${\cal K}$ is a holomorphic vector field on ${\cal PT}$ 
corresponding to a conformal Killing vector $K$.

Below we shall state the Penrose result extended to the Einstein and 
hyper-Hermitian cases:
\begin{prop}
\label{twmaci}
Let $\cal PT$ be a three-dimensional complex manifold with a
four-dimensional family of rational curves (invariant under a complex
conjugation with fixed points) with
normal bundle ${\cal O}(1)\oplus {\cal O}(1)$.
Then the moduli space ${\cal M}$ of these sections is equipped with 
an ASD conformal structure $[g]$ of signature $(++--)$. Conversely given an ASD
four-manifold there will always exists a corresponding twistor space.
Moreover ${\cal M }$ is:
\begin{itemize}
\item Hyper-K\"ahler iff  there exists 
a projection $\mu :{\cal PT}\longrightarrow \CP^1$, and  each fibre
of this projection is equipped with an $\mu^*\O(2)$ valued symplectic form
{\em{\cite{Pe76}}} (equivalently, we can require that the canonical
bundle $\kappa$ of ${\cal PT}$ is $\kappa=\mu^*{\cal O}(-4)$).
\item Hyper-Hermitian iff  there is 
a projection $\mu :{\cal PT}\longrightarrow \CP^1$ {\em{\cite{D98}}}.
\item Einstein ($R_{ab}=\Lambda g_{ab}$) iff there exists
a contact structure $\tau\in\Lambda^2(T^*{\cal PT})\otimes \O(2)$,
where now $\O(2)=\kappa^{-1/2}$, and $\kappa$ 
is the canonical bundle $\Omega^{3}$, 
such that $\tau\wedge\d\tau=\Lambda\xi$ where $\xi \in
\Omega^{3}\otimes \kappa^{-1}$ {\em{\cite{W80}}}.
\end{itemize}
\end{prop}
\subsubsection{Construction of the two--form}
\label{construction_form}
Consider an ASD four--manifold $({\cal M}, [g])$.  Define
the non-projective twistor space, ${\cal T}$, to be the total space of
the line bundle $\kappa^{1/4}\rightarrow {\cal PT}$ where 
$\kappa=\Omega^3$ is the canonical bundle.  
In the conformally-flat case ${\cal T}$ is the tautological
line bundle ${\cal O}(-1)$, i.e.\ $\C^4\mapsto\CP^3$, and we will also
use this notation, ${\cal T}={\cal O}(-1)$ in the curved case.  The
nonprojective spin bundle $S_{A'}\mapsto {\cal M}$ is defined to be the
total space of the pullback of this line bundle to the correspondence
space ${\cal F}={\cal M}\times\CP^1$.  Clearly $S_{A'}={\cal M}\times\C^2$.  
The fibration $q:S^{A'}\mapsto {\cal T}$ is spanned by a lift
of the twistor distribution or Lax pair.  The non-projective spin
bundle is the total space of a line bundle, which we will also denote
by ${\cal O}(-1)$, over ${\cal F}$.  (Note that in the hyper-Hermitian
case the line bundles $\O(n)$ just defined will {\em not} be the same as
$\mu^* \O(n)$ unless $({\cal M}, [g])$ is in fact hyper-Kahler.)

The space ${\cal T}$ admits an Euler vector field $\Upsilon$ being the
total space a of line bundle, and a tautological three-form, $\xi$ the
pullback of the tautological three-form on $\kappa$.  These satisfy
${\cal L}_\Upsilon\xi=4\xi$.  Let $\phi=\d \xi$, then
$\xi=4\phi(\Upsilon, ..., ...)$. $\xi$ can be thought of as a form on
${\cal PT}$ with values in the dual canonical bundle $\kappa^*$.

We now impose a symmetry: let $K, \Kt$, and ${\cal K}$ be
respectively: a conformal Killing vector on ${\cal M}$, its lift to
the correspondence space ${\cal M}\times\CP^1$, and the holomorphic
vector field on ${\cal T}$ which is the push-forward of $\Kt$.
\begin{prop}
The two form $\Smt:=q^*\phi({\cal K}, \Upsilon, ..., ...)\in \Lambda^2(T^*S^{A'})$
 satisfies
\be
\label{closureetc}
\Smt\wedge\Smt=0,\;\;\;\;\;\d\Smt=\beta\wedge\Smt\;\;\;\;\;{\cal
  L}_{\Kt}\Smt=0 \ee for some one-form $\beta$ homogeneous of degree
$0$ in $\pi^{A'}$.
\end{prop}
{\bf Proof:} 
It follows from the definition of $\Smt$  that the
integrable twistor distribution  belongs is the kernel of 
$\Smt$. Therefore equations (\ref{closureetc}) follow from 
Frobenius' theorem.
The one-form $\beta$ is defined up to  the addition of $\d(\ln{\sigma})$
where $\sigma$ is a twistor function homogeneous of degree 0. 
\koniec
\smallskip

From $ {\cal L}_{\Upsilon}\Smt=4\Smt$ and ${\Upsilon}\hook\Smt=0 $ it
follows that $\Smt$ descends to ${\cal F}$ where it takes values in
$\O(4)$. Note however that $\d \Smt$ does not descend as
${\Upsilon}\hook\d\Smt={\cal L}_{\Upsilon}\Smt\neq 0$.  To
differentiate $\Smt$ on ${\cal F}$ we need a nonzero section of
$\O(4)$ in order to dehomogenise $\Smt$.  When $({\cal M}, g)$ is ASD
Einstein or vacuum we can find a section of ${\cal O}(4)$ to
dehomogenise $\Smt$. This section necessarily has zeroes, and so
equivalently, this requires the existence of a divisor description of
the dual canonical bundle.  This can be seen from the twistor construction.
\begin{itemize}
\item
{\bf Vacuum case:} The twistor space fibres over $\CP^1$ and so we can
pull back $\pi\cdot\d\pi$ to ${\cal PT}$.
\begin{figure}
\caption{Divisor on a mini-twistor space.}
\label{devisor}
\begin{center}
\includegraphics[width=8cm,height=11cm,angle=0]{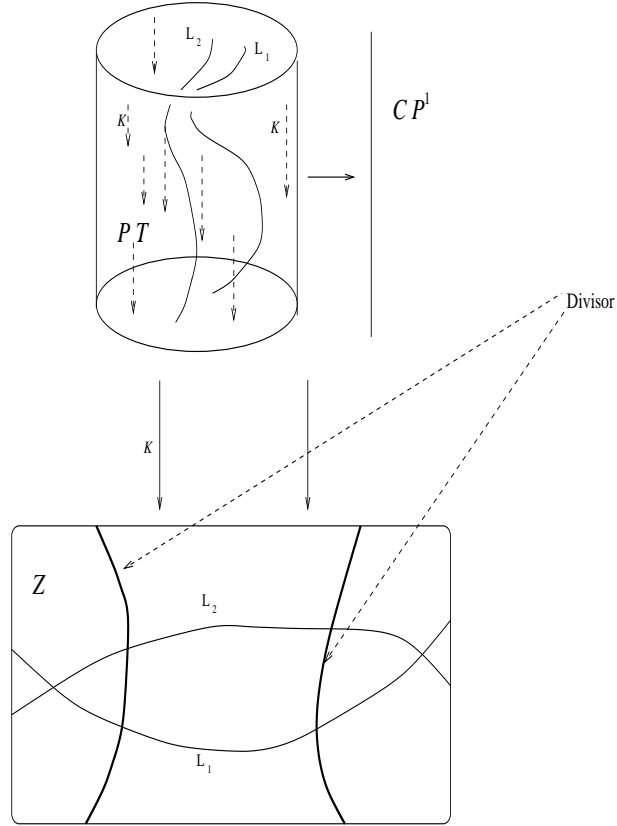}
\end{center}
\end{figure}
Let ${\cal K}$ be a holomorphic vector field on ${\cal PT}$ such that
${\cal L}_{\cal K}\Sm_{\l}=\eta\Sm_{\l}$ (${\cal K}$
corresponds to a Homothetic Killing vector on ${\cal M}$).  The
function $D:={\cal K}\hook\pi\cdot\d\pi$ is a section of $\O(2)$ and
the two-form $D^{-2}{\cal K}\hook\xi$ descends to the mini-twistor
space ${\cal Z}$.
\item
{\bf Einstein case:}
Let ${\cal PT}_E$ be the projective 
twistor space corresponding to a solution of the ASD Einstein equations.
It is equipped with a
contact structure $\tau\in\Lambda^2(T^*{\cal PT}_E)\otimes \O(2)$
such that $\tau\wedge\d\tau=\Lambda\xi$.
$\d \tau$ defines a holomorphic symplectic structure on the
non-projective twistor space ${\cal T}_E$.  If $K$ is a Killing vector
on an ASD Einstein manifold then the corresponding holomorphic vector
field on the non-projective twistor space is Hamiltonian with respect
to $\d \tau$.  To see this, define a section of $O(2)$ by $D:={\cal
K}\hook\tau$.  We have $ \d D={\cal L}_{\cal K}\tau-{\cal K}\hook
\d\tau=-{\cal K}\hook \d\tau$ as ${\cal K}$ is a symmetry.
\end{itemize}
On the projective spin bundle ${\cal F}$ define
\[
\Pi:=D^{-2}\Smt.
\]
We have the following result:
\begin{prop}
\label{twoform2}
The two-form $\Pi$ is well defined on the Einstein--Weyl correspondence
space ${\cal F}_W$. It satisfies  
\be
\label{EWform}
\d \Pi=0,\;\;\;\;\;\;\Pi\wedge \Pi=0,
\ee
where $\d=\d x^i\otimes\p_i+\d \lt\otimes \p_{\lt}$ 
is the exterior derivative on  ${\cal F}_W$. 
Any two linearly independent vectors $L_{A'}$ such that $L_{A'}\hook
S=0$ form a Lax pair for the EW equations.
\end{prop}
{\bf Proof.} The simplicity follows from $\Smt\wedge\Smt=0$.  In the
vacuum case the two form 
\be
\label{pull_back_form}
{\Pi}=q^*\frac{{\cal K}\hook \xi}{{\cal K}\hook(\pi\cdot\d\pi)}
\ee
is a pull back of a closed and simple form on ${\cal PT}$.
In the Einstein case 
\[
\Pi= D^{-2}q^*{\cal K}\hook(\Lambda\tau \wedge \d\tau)=\d(\Lambda\tau/D).
\]
Therefore Einstein--Weyl metrics which come from ASD Einstein and
hyper-K\"ahler four manifolds give rise to the same structure on the
reduced spin bundle.
The form $\Pi$ descends to ${\cal F}_W$
because $\Kt\hook\d \Pi=0$ and $\d(\Kt\hook \Pi)=0$.
\koniec 
\smallskip

\noindent
 {\bf Remark.}
In \cite{Ta90} certain dispersionless integrable systems were
expressed in terms of $\Pi$ satisfying (\ref{EWform}).
\smallskip

The two form $\Smt$ can be equivalently constructed from the data on
${\cal M}$ as follows.  Let $K$ be a Killing vector on a general ASD
conformal manifold $({\cal M}, [g])$, and let $\Xi$ be a volume form
on the non-projective primed spin bundle $S^{A'}$.  Define the two
form on $S^{A'}$ \be
\label{twoform1}
\Smt:=\Xi(L_0, L_1, \Kt, \Upsilon_{\Xi}, ..., ... ).
\ee
Here $\Upsilon_{\Xi}=\pi^{A'}/\p \pi^{A'}$ is the Euler vector field on
$S^{A'}$, $L_A$ is the twistor distribution,
and $\Kt$ is a Lie lift of $K$ to $S^{A'}$. 
Now assume that $({\cal M}, g)$ is also vacuum. Consequently
$\nabla_{AA'}{K^A}_{B'}=const$ and the spin bundle is equipped with a canonical
divisor\footnote{We assume that $\nabla_{AA'}{K^A}_{B'}\neq 0$. 
If $\nabla_{AA'}{K^A}_{B'}=0$ then $K$ is triholomorphic and a
section of $\O(2)$ which descends to the reduced spin bundle
is $(\iota\cdot\pi)^2$ where $\iota_{A'}$ is any constant spinor.}
$D:=\pi^{A'}\pi^{B'}\nabla_{AA'}{K^A}_{B'}\in\O(2)$ which descends to
the reduced spin bundle\footnote{By the reduced spin bundle
(correspondence space) we
mean the space of orbits of $\Kt$ in $S^{A'}$ (in ${\cal F}$).}
(Figure \ref{devisor}).
It is easy to prove  that now
\begin{eqnarray}
\label{beta}
\Smt&=&\pi_{A'}\pi_{B'}\pi_{C'}\pi_{D'}\phi^{A'B'}\Sm^{C'D'}
+\pi_{A'}\pi_{B'}\pi_{C'}\d\pi^{C'}\wedge
(K\hook\Sm^{A'B'}),\nonumber\\
\beta&=&\frac{4\phi_{A'B'}\pi^{A'}\d \pi^{B'}}{\pi_{A'}\pi_{B'}\phi^{A'B'}}=
\d\ln{D^2}\nonumber\\
\Pi&=&\d\l\wedge\frac{K\hook\Sm(\l)}{D^2}-\frac{\Sm(\l)}{D},\qquad\mbox{where}\qquad
\Sm(\l)=\pi_{A'}\pi_{B'}\Sm^{A'B'}.
\end{eqnarray}  
From the last formula it follows that to construct $\Pi$ one should rewrite $\Sm(\l)/D$
in the coordinates in which $K=\p_t$, and then replace all $\d t$s by the differentials of 
a suitably defined invariant spectral parameter.

{\bf Example.} We shall now illustrate the construction of $\Pi$ with
a simple example. 
Let $2\d w\d \tw-2\d z\d \tz$ be a flat metric on $\R^{2,2}$ and let
$K=z\p_z-\tz\p_{\tz}$ be a Killing vector.
The flat twistor distribution and the lifted symmetry are:
\[
L_0=\p_{\tw}-\l\p_z,\;\;\;L_1=\p_{\tz}-\l\p_w,\;\;
\Kt=z\p_z -\tz\p_{\tz}+\l\p_{\l}.
\]
The volume form on ${\cal F}$ and the two-form $\Sm(\l)$
are given by
\[
\Xi=\d\l\wedge \d z\wedge \d\tz \wedge \d w \wedge \d \tw, \;\;
\Sm(\l)=-\l^2\d \tw \wedge \d\tz +\l( \d w \wedge \d \tw -\d z\wedge \d\tz)
+\d w\wedge \d z.
\]
In the covariantly constant frame we introduce
$
2r:=\mbox{ln}(z\tz),\;2\phi:=\mbox{ln}(z/\tz)$, so that 
$\Kt=\p_{\phi}+\l\p_{\l}$.
In these coordinates
\[
\Sm(\l)=-\l^2e^{r-\phi}\d\tw\wedge(\d r-\d\phi)+\l(\d w\wedge
\d\tw+2e^{2r}\d r\wedge \d\phi)+e^{r+\phi}\d w\wedge(\d r+\d\phi)
\]
and (from (\ref{beta}))
\be
\label{flattwoform}
\Pi=e^r(\d\tw\wedge \d\lt +\lt^{-2}\d w\wedge \d\lt+\lt \d\tw\wedge \d r
-{\lt}^{-1}\d w\wedge \d r)
+2{\lt}^{-1}e^{2r}\d r\wedge \d\lt -\d w\wedge \d\tw
\ee
where $\lt=\l e^{-\phi}$ is an invariant spectral parameter.

 The two form $\Pi$ can be also obtained as a 
pull-back from ${\cal PT}$.
 Local inhomogeneous coordinates on ${\cal PT}$ pulled back to ${\cal F}$  are given by
$
(\l,\;\mu^1=\l\tw+z,\;\mu^0=\l\tz+w).
$
The holomorphic vector field on  ${\cal PT}$ is
${\cal K}=\mu^0\p_{\mu^0}+\l\p_{\l}$. From (\ref{pull_back_form}) we have
\[
q^*({\cal K}\hook(\d\l\wedge \d\mu^0 \wedge \d\mu^1)=
(\mu^0 \d\l-\l \d\mu^1)\wedge
\d\mu^1 =\l^2\d\mu^1\wedge \d(\mu^0/\l).
\]
Thus
\[
\Pi=\d\mu^1\wedge \d(\mu^0/\l)=\d P\wedge\d Q
\]
which agrees with (\ref{flattwoform}).
Here ${P}=\tw+{\lt}^{-1}e^r$ and ${Q}=\lt e^r+w$ are coordinates on
mini-twistor space pulled back to the reduced spin bundle.

\section{Twistor theory of the dKP Einstein-Weyl structures}
\label{twistor_theory}

Here we give an account of the twistor theory of the dKP EW metrics,
and the dKP equation
(some connections between a twistor theory and the  
dKP equations have  been discussed in \cite{GT}).
We shall also characterise all four dimensional 
hyper-K\"ahler and ASD Einstein metrics that give rise to 
the dKP EW structures.

Define
the non-projective twistor space, ${\cal Y}$ corresponding to a Weyl
space ${\cal W}$, 
to be the total space of
the line bundle $\kappa^{1/4}\rightarrow {\cal Z}$ where 
$\kappa=\Omega^2$ is the canonical bundle of ${\cal Z}$.  
The nonprojective spin bundle $S_{A'}\mapsto {\cal W}$ is the rank two
vector bundle defined to be the
total space of the pullback of this line bundle to the correspondence
space ${\cal W}\times\CP^1$.  The fibration $q:S^{A'}\mapsto {\cal Y}$
is spanned by a lift of 
the mini-twistor distribution $L_{A'}$ (\ref{LAX_PAIR}).

Any shear-free null geodesic congruence of the Einstein-Weyl structure
determines a one-dimen\-sional sub-manifold in $\cal Z$
(this is a reduction of the 4-dimen\-sional Kerr theorem).
A codimension--one submanifold determines a line bundle $[D]$ by
the divisor construction; $[D]$ admits a section $D$ that vanishes
precisely on the given submanifold.

When the Einstein--Weyl geometry arises from a solution of the dKP
equation the dual canonical bundle $\kappa^{-1}$ of the minitwistor space
admits a fourth root that is given by the divisor construction, that
is it admits a section $D$ that vanishes on a codimension-one subset.
In general, as seen above, if the Einstein-Weyl geometry is a
reduction of an ASD Einstein, or hyper-K\"ahler four-manfiold, 
then $\kappa^{-1/2}$ admits a
section whose zero set will generally have two components in the
neighbourhood of a line.  For an Einstein-Weyl dKP solution, the two
`divisor curves' in Fig (\ref{devisor}) degenerate to one curve.
This observation gives rise to a twistor characterisation of solutions
to the dKP equation

\begin{prop}
\label{dKPtwistor}
There is a one to one correspondence 
between Einstein-Weyl spaces obtained from solutions to the dKP equation
and two-dimensional complex manifolds with
\begin{itemize}
\item A three parameter family of rational curves with normal bundle
$\O(2)$.
\item
A global section $l$ of $\K^{-1/4}$, where $\K$ is the canonical bundle.
\end{itemize}
In order to obtain a real Einstein-Weyl structure, we require an
antiholomorphic involution fixing a real slice, leaving a rational
curve invariant and leaving the section of
$\K^{-1/4}$ above invariant.
\end{prop}
{\bf Proof.}  
The global section $l$ of $\kappa^{-1/4}$, when pulled back to $S_{A'}$
determines a homogeneity degree one function on each fibre of $S_{A'}$
and so must, by globality, be given by $l=\iota^{A'}\pi_{A'}$ and since
$l$ is pulled back from twistor space, it must satisfy $L_{A'}l=0$.
This implies $\widetilde{D}_{A'(B'}\iota_{C')}=0$, and 
(after some algebraic manipulations)
\[
\widetilde{D}_{A'B'}\iota^{C'}=0,
\]
where $\widetilde{D}$ is a covariant weighted derivative.

Therefore the null vector field $l^{a}=\iota^{A'}\iota^{B'}$ is covariantly
constant. The Lemma \ref{lemma_weight} implies that the conformal 
weight of $\iota^{A'}$ is $-1/4$ and hence that of
$l^{a}$ is $-1/2$.  This weight can be deduced from the correspondence
as follows: the two form ${\widetilde\Sm}=\pi_{A'}\pi_{B'}e^{A'B'}\wedge
\varepsilon^{C'D'}\pi_{C'}\d\pi_{D'}$ has conformal weight $0$ on
$S^{A'}$. $e^{A'B'}$ has weight 0, and $\varepsilon^{A'B'}$ weight
$-1$ so $\pi_{A'}$ has weight $1/4$. The global section
$\pi_{A'}\iota^{A'}$ is weightless so 
the weight of $\iota^{A'}$ is $-1/4$. 
Hence by Proposition \ref{todth} the corresponding Einstein-Weyl space
arises from a solution to the dKP equation.

Conversely, given a solution to (\ref{dKP}) one can obtain ${\cal Z}$
as a factor space of ${\cal W}\times \CP^1$ by the distribution
(\ref{dKPlax}) and the covariant constant weighted null vector
$l^a=\iota^{A'}\iota^{B'}$ gives
rise to the section $l=\iota^{A'}\pi_{A'}$ of $\K^{-1/4}$ 
\koniec \smallskip

\noindent
{\bf Remark:} Note that there is not a $1-1$ correspondence between
such twistor spaces and solutions to the dKP equation on account of
the coordinate freedom (\ref{coordtrans}) and (\ref{conftrans}).  The
coordinate choices implicit in a solution to the dKP equation can be
encoded on the twistor space in the choice of the coordinates near the
divisor as follows.

Let $\hat{P}, \hat{Q}$ be local coordinates on a
neighbourhood of the divisor in $\cal Z$ such that $\hat{Q}=0$ on the
divisor and, setting $Q=\hat{Q}^{-1}, P=\hat{P}/\hat{Q}^2$ on the
complement of the divisor,we have
\[
\Pi = \d P\wedge\d Q=-\hat{Q}^{-4}\d \hat{P}\wedge\d \hat{Q}.
\]
Consider a graph of a rational curve ${\hat P}({\hat Q})$. Parametrise
the curve by $(t, y, x)$ as follows:
\[
t:=\hat{P}|_{\hat{Q}=0},\qquad 
y:=\frac{\d \hat{P}}{\d \hat{Q}}|_{\hat{Q}=0},\qquad
x:=\frac{1}{2}\frac{\d^2 \hat{P}}{\d \hat{Q}^2}|_{\hat{Q}=0}.
\]
Therefore the local 
coordinates $P, Q$ have the following expansion near $\lt=\infty$ 
\[
Q:=\lt +\sum_{i=1}^{\infty}u_i\lt^{-i},\qquad P=\sum_{i=1}^{\infty}w_iQ^{-i}
+x+Qy+Q^2t
\]
(after performing  an $SL(2, \C)$ transformation and choosing a spin
frame such that the constant term in the Laurent expansion of $Q$
vanishes).
When we pull the mini-twistor coordinates back to ${\cal F}$, then 
$u_i, w_i$ become functions of $(x, y, t)$. 
The functions $P$ and $Q$ are solutions of Lax equations $L_{A'}P=L_{A'}Q=0$.
They form a local Darboux atlas as $\Pi=\d P\wedge\d Q $, where $\Pi$ is given
by (\ref{dKPtwoform}).
\[
\Pi=\d x\wedge\d \lt+\d y\wedge
\d(\frac{{\tilde\lambda}^2}{2}+u_1) +
\d t\wedge \d (\frac{{\tilde\lambda}^3}{3}+
\tilde{\lambda}{u_1}+w_1).
\]
The poles of $\Pi$ occur on the divisor.
Now $\Pi$ is a pull back of a two-form from a two-dimensional
manifold. Therefore is satisfies $\Pi\wedge\Pi=0$, which yields
${w_1}_x={u_1}_y$ and
the dKP equation (\ref{dKP}) for $u_1$. 

Thus, a solution to the dKP equation corresponds to a EW mini-twistor space
as described in Proposition \ref{dKPtwistor} together with  a Darboux
coordinate system as above on the third formal neighbourhood of the divisor.
[It seems likely that the Benney hierarchy will similarly correspond
to the EW dKP minitwistor space as above together with the Darboux
coordinate system on a neighbourhood of the divisor defined now to all orders.]

\smallskip

Now we are in a position to give a characterisation of the hyper-K\"ahler
metrics (\ref{HCdKP}).
\begin{prop}
Let $g$ be an indefinite  hyper-K\"ahler metric with a symmetry $K$ satisfying 
$\d K_+\wedge \d K_+=0$. Then $g$ is locally of the form 
{\em(\ref{HCdKP}\em)}.
\end{prop}
{\bf Proof.}
Let ${\cal K}$ be a vector field (corresponding to $K$) 
on a  twistor space of $({\cal M}, g)$. The divisor
\[
{\cal K}\hook \pi\cdot\d\pi=\pi_{A'}\pi_{B'}\phi^{A'B'}
\] descends to
the minitwistor space. If $\d K_+$ is null then 
$\phi_{A'B'}=(1/2)\nabla_{AA'}K^{A}_{B'}=\iota_{A'}\iota_{B'}$ for some
constant spinor $\iota^{A'}$. Therefore $\pi\cdot \iota$ on ${\cal PT}$
defines a divisor in $\cal Z$. It takes values in $\K^{-1/4}$ because
the canonical bundle of ${\cal PT}$ is the square of the pullback of
the canonical
bundle of $\CP^1$. The assumptions of Proposition
\ref{dKPtwistor} are satisfied and so the EW structure corresponding
to ${\cal Z}$ is of the form (\ref{EWdkp}).
Therefore it follows from Proposition \ref{prop_JT} 
that the metric $g$ is given by
\[
g=\Omega(\t{V}(\d \t{y}^2-4\d \t{x}\d \t{t}-4\t{u}\d
\t{t}^2)-\t{V}^{-1}
(\d \t{z}+\t{\a})^2)=\Omega \tilde{g},
\]
where $\t{u}(\t{x}, \t{y}, \t{t})$ is a solution to dKP 
$(\t{V}, \t{\a})$ is a solution to the monopole equation
(\ref{EWmonopole}), and $\Omega$ is a conformal factor.
Calculating the scalar curvature of the metric $\tilde{g}$ yields 
\[
\tilde{R}=8(\t{V}_{\t{y}\t{y}}-\t{V}_{\t{x}\t{t}}+(\t{u}\t{V})_{\t{x}\t{x}})\t{V},\]
and so $\tilde{R}=0$ because $\t{V}$ satisfies (\ref{lindKP}). However
the  metric $g$ is hyper-K\"ahler, therefore its scalar curvature also
vanishes. As a consequence  we  deduce that $\Omega=\Omega(\t{t})$.
Now we can use the coordinate freedom (\ref{conftrans}) to absorb
$\Omega$ in the solution to the dKP equation. This yields
\be
\label{scalarflat}
g=(V(\d y^2-4\d x\d t-4u\d t^2)-V^{-1}(\d z+\a)^2)=\Omega \tilde{g},
\ee 
where $({V}, {\a})$ is another solution to the monopole equation.
In section \ref{HKsection} we showed that  
this metric is hyper-K\"ahler metric if $V$ is a multiple of $u_x$.

 Consider the metric (\ref{scalarflat}) with an arbitrary
monopole $V$ (an arbitrary solution to the linearised dKP equation
\ref{lindKP}). The self-dual derivative of the isometry $K=\p_z$
is given by $\phi_{A'B'}=(u_x/V)\iota_{A'}\iota_{B'}$, for some
constant spinor 
$\iota_{A'}$. The well known identity  
$
\nabla_a\nabla_bK_c=R_{bcad}K^d
$
and the vacuum condition yield $\nabla_{a}\phi_{B'C'}=0$. Therefore
(\ref{scalarflat}) is  hyper-K\"ahler iff $u_x/V=const$.
\koniec 

\smallskip
\noindent
{\bf Remarks:}
\begin{itemize}
\item This Proposition corrects an omission made
in the classification \cite{FP79} of complexified hyper-K{\"a}hler
spaces with symmetry. In the Appendix we shall demonstrate explicitly
that the dKP equation is a reduction of the second heavenly equation 
considered in  \cite{FP79}. 
\item  Metrics (\ref{scalarflat}) with $V\neq const\times u_x$ are
not vacuum, but they 
admit a covariantly constant real spinor. The full characterisation
of these metrics will be given in our subsequent paper.
\end{itemize}
\begin{prop}
All EW structures which arise from indefinite
ASD  Einstein metric with a symmetry $K$ satisfying 
$\d K_+\wedge \d K_+=0$ are locally of the form 
{\em(\ref{EWdkp}\em)}.
\end{prop}
{\bf Proof.}
The canonical divisor  $D:={\cal K}\hook\tau$ (where $\tau$ is the
contact structure) descends to a mini-twistor space. Because 
$\d K_+$ is null the square root of $D$
exists and takes its values in $\K^{-1/4}$.
\koniec

\section{Symmetry reductions of hyper-K\"ahler metrics in $2+2$ signature}
Symmetry reductions of the hyper-K\"ahler condition on a  real 
four-dimensional Riemannian metric have been completely classified: 
\begin{itemize}
\item If the symmetry is tri-holomorphic, 
then the corresponding metric belongs to the Gibbons--Hawking class
\cite{GH78}, and is given by a solution to the Laplace equation 
in three dimensions. The resulting Einstein--Weyl structures are
trivial, and their mini-twistor space is $T\CP^1$.
\item Hyper-K\"ahler metrics with non-triholomorphic Killing
vectors are given by solutions to the $SU(\infty)$ Toda equation \cite{FP79}.
The corresponding EW structures \cite{W90} 
are characterised by the existence of
a shear-free, twist-free geodesic congruence \cite{T95}. 
Mini-twistor spaces are in this case equiped with a  
canonical divisor (two one-dimensional complex sub-manifolds)
taking its values in $\O(2)$ \cite{L91}. In 
\cite{C00} EW Toda structures were characterised in terms of weighted
vector fields.

\item Hyper-K\"ahler metrics with tri-holomorphic conformal symmetries
yield a class of EW structures (called hyper-CR EW structures)
characterised by the existence of a sphere of shear-free, 
divergence-free geodesic congruences \cite{GT98}. 
The corresponding mini-twistor spaces are fibred over $\CP^1$.
\item Hyper-K\"ahler metrics with non-tri-holomorphic, conformal
symmetry (and the resulting EW structures)
are given by solutions to a certain second order integrable equation
in three dimensions \cite{DT99}. This equation gives
$SU(\infty)$-Toda and hyper-CR Einstein-Weyl structures as limiting
cases. The EW structures arising from conformal, non-tri-holomorphic
reductions are characterised by the existence of a shear-free
geodesic congruence for which the twist is a constant 
multiple of the divergence \cite{CP99}.
\end{itemize}
The above list is not complete if one considers   
Hyper-K\"ahler metrics in $(++--)$ signature. 
The existence of null structures of various kinds allows 
two additional types of symmetries:
\begin{itemize}
\item
Hyper-K\"ahler metrics  
for which the self-dual part of a derivative of
a Killing vector is null correspond to solutions of the dispersionless
Kadomtsev--Petviashvili equation (\ref{dKP}). 
The corresponding EW structures are characterised by the existence of
a constant weighted vector. The minitwistor spaces are such that 
the line bundle 
$\K^{-1/4}$ admits a section, where $\K$ is the canonical line bundle.
The above statements have been proved in this paper.
\item
Hyper-K\"ahler metrics with conformal Killing vectors 
for which the self-dual part of a derivative of
a conformal Killing vector is null. 
\end{itemize}
The last possibility has not yet been investigated. The EW
spaces will be given by a generalisation of the dKP equation.
We intend to study this generalisation, and the
corresponding EW geometries in a subsequent paper.

\section{Outlook: a twistor theory for the full KP equation?}
A combination of the dispersive limit of dKP 
with the twistor picture 
suggests
a candidate for a twistor space for the full KP equation (\ref{KP})
(cf the similar proposal in \cite{S95}).

Let $x$ be a coordinate on a configuration space $Q$,
and let $\lt$ be the corresponding momentum. 
The extended six-dimensional 
phase-space $T^*(Q\times \R^2)$ is coordinatised by 
$x^i=(x, y, t), p_i=(\lt, H_2, H_3)$.
Restrict the symplectic form $\Pi$ 
on $T^*(Q\times \R^2)$ to the four-dimensional correspondence space
${\cal F}^4$ 
obtained by putting $H_r:=H_r(x^i, \lt)$, $r=2,3$.
The (complexified) space ${\cal F}^4$ is foliated by sub-manifolds
whose tangent vectors annihilate the symplectic form,
which gives rise to a projection $p:{\cal F}\longrightarrow {\cal Z}$
such that $\Pi$ descends to a symplectic form on ${\cal Z}$.
The two-dimensional complex manifold ${\cal Z}$ is the mini-twistor space for the extended
configuration space $Q\times \R^2$ with its dKP Einstein--Weyl structure.
It is believed that the Moyal quantisation of $T^*(Q\times \R^2)$
gives rise to  
the full KP equation. This suggests the conjecture that there 
exists a correspondence between solutions to the full KP equation
and the Moyal deformations  of 
${\cal Z}$. 

It will be instructive to compare this approach to the twistor
constructions for the full KP equations described in 
\cite{M95}, and \S12.6 of \cite{MW96}.

\section{Acknowledgements}
We thank David Calderbank for valuable comments which resulted
in many improvements.  
Maciej Dunajski
would like to thank {\em Centro de Investigacion y de Estudios
Avanzados} in Mexico, where part of this work was done, for its
financial support (3697 E, Proyecto de CONACYT).  MD is also grateful
to Bogdan Mielnik, Maciej Przanowski and Jerzy Pleba\'nski for their
warm hospitality.  LJM would like to acknowledge support from NATO
collaborative Research Grant number CRG 950300.

Some results in sections 3 and 5 were obtained during the workshop
{\em Spaces of geodesics and complex methods in general relativity and
geometry } held in the summer of 1999 at the Erwin Schr{\"o}dinger
Institute in Vienna.  We wish to thank ESI for the hospitality and for
financial assistance.

\section{Appendix}
Here we shall demonstrate (by an explicit calculation) that the 
dKP equation (\ref{dKP}) is a reduction of the second heavenly equation
by  a Killing vector with a null self-dual derivative.

 Let $\Th(z, t, q, y)$ satisfy \cite{Pl75}.
\be
\label{secondeq}
\Th_{zy}-\Th_{tq}+\Th_{qq}\Th_{yy}-\Th_{qy}^2=0.
\ee
Then 
\be
\label{Plebanm}
g=2(\d z\d y +\d q\d t -\Th_{qq}\d z^2-\Th_{yy}\d t^2+2\Th_{yq}\d z\d t)
\ee
is a hyper-K\"ahler metric. All hyper-K\"ahler metrics can locally
be put in the form (\ref{Plebanm}).

Let $K$ be a Killing vector such that $\d K_+\wedge\d K_+=0$.
There is no loss of generality \cite{FP79} in choosing $K=\p_z-2z\p_q$,
in which case $\d K_+=2\d t\wedge\d z$.

The Killing equations yield $({\cal L}_K\Th)_{yy}=
({\cal L}_K\Th)_{qq}=0, ({\cal L}_K\Th)_{yq}=1$. 
They integrate to
\be
\label{dKPtheta}
\Th=zqy+yA(z, t)+qB(z, t)+C(z, t) +G(y, t, q+z^2).
\ee
The function $C$ is pure gauge and can be set to zero without loss of
generality.
Imposing (\ref{secondeq}) gives two equations:
the first is $A_z+B_t=2z^2$, and we can deduce, without loss of
generality, that $A=z^3, B=-z^2t$,
and the second is
\be
\label{dKPlegandre}  
-u-G_{tu}+G_{yy}G_{uu}-G_{yu}^2=0,\qquad\mbox{where}\qquad u=-(q+z^2).
\ee
The last equation is equivalent to the dKP equation. To see this 
write (\ref{dKPlegandre}) as a closed system 
\begin{eqnarray}
\d G&=&G_u\d u+G_t\d t+G_y\d y,\nonumber\\
\label{legen1}
0&=&-u\d y\wedge\d t\wedge \d u +\d G_u\wedge\d y\wedge \d u
-\d G_y\wedge\d G_u\wedge \d t.
\end{eqnarray}
Now rewrite the first equation as
$
\d(G-uG_u)=G_t\d t+G_y\d y-u\d G_u,
$
and perform a Legendre transform  
\[
x:=G_u,\qquad u=u(t, y, x),\qquad
H(t, y, x):=-G(t, y, u(t, y, x))+xu(t, y, x).
\]
The relation 
$
\d H=H_t\d t+H_x\d x+H_y\d y
$
implies $H_t=-G_t, H_y=-G_y, H_x=u$.
Equation (\ref{legen1}) yields
\[
-H_x\d y\wedge\d t\wedge\d H_x+\d x\wedge\d y\wedge\d H_x+\d H_y\wedge\d x
\wedge\d t=0,
\]
which is equivalent to 
\be
\label{almostdKP}
H_xH_{xx}-H_{xt}+H_{yy}=0.
\ee
Taking the $x$ derivative of the above equation and using
$H_x=u$ yields
\[
u_{xt}-uu_{xx}-u_{x}^2=u_{yy}
\]
which is the dKP equation.
To calculate the metric differentiate the relation $x=G_u$ 
with respect to $x$ and $H_y=-G_y$ with respect to $y$,
\[ 
1=G_{uu}u_x,\qquad 0=G_{uy}+G_{uu}u_y,\qquad 0=G_{ut}+G_{uu}u_t,\qquad
G_{yy}=\frac{{u_y}^2}{u_x}+uu_x-u_t
\]
(we also used (\ref{almostdKP})).
Therefore (from (\ref{dKPtheta})) we have
\[
\Th_{yy}=\frac{{u_y}^2}{u_x}+uu_x-u_t,\qquad
\Th_{yq}=\frac{u_y}{u_x}+z,\qquad\Th_{qq}=\frac{1}{u_x}.
\]
The metric (\ref{Plebanm}) in terms of $u(x, y, t)$ is
\begin{eqnarray*}
g&=&2(-u_x\d x\d t+\d z\d y+2\frac{u_y}{u_x}\d z\d t-u_y\d y\d t
-(uu_x+\frac{u_y^2}{u_x})\d t^2-\frac{1}{u_x}\d z^2)\\
 &=&\frac{u_x}{2}(\d y^2-4\d x\d t-4u\d t^2)
-\frac{2}{u_x}(\d z-\frac{u_x\d y}{2}-u_y\d t)^2
\end{eqnarray*}
which is (\ref{HCdKP}).


\end{document}